\documentclass[12pt]{article}
\usepackage{amssymb,amsfonts,amsmath,amsthm}
\usepackage[dvips]{graphicx}
\newcommand{\ol}{\overline}
\newcommand{\one}[1]{\mbox {\bf 1}_{\{#1\}}}
\newcommand{\odin}{\mbox {\bf 1}}
\newcommand{\pp}{{\mathbb P}}

\newcommand{\ra}{\rightarrow}

\newcommand{\witi}{\widetilde}

\newcommand{\fracd}[2]{\frac {\displaystyle #1}{\displaystyle #2 }}

\newcommand{\suml}{\sum\limits}

\newcommand{\prodl}{\prod\limits}
\newcommand{\zz}{{\mathbb Z}}
\newcommand{\nn}{{\mathbb N}}
\newcommand{\ee}{{\mathbb E}}
\newcommand{\rr}{{\mathbb R}}
\newcommand{\cals}{{\mathcal S}}
\newcommand{\call}{{\mathcal L}}
\newcommand{\calt}{{\mathcal T}}

\newcommand{\calf}{{\mathcal F}}
\newcommand{\calg}{{\mathcal G}}
\newcommand{\calb}{{\mathcal B}}

\newcommand{\veps}{\varepsilon}
\newcommand{\beq}{\begin{eqnarray*}}
\newcommand{\feq}{\end{eqnarray*}}
\newcommand{\beqn}{\begin{eqnarray}}
\newcommand{\feqn}{\end{eqnarray}}
\newcommand{\be}{\begin{equation*}}
\newcommand{\fe}{\end{equation*}}
\newcommand{\ben}{\begin{equation}}
\newcommand{\fen}{\end{equation}}

\newcommand{\as}{\mbox{a.s.}}
\newcommand{\f}{{\mathcal F}}

\newcommand{\proda}{{\mbox {\Large $\Pi$}}}
\newcommand{\citethm}[2]{[{#1\if@tempswa,#2\fi}]}
\makeatletter \@addtoreset{equation}{section} \makeatother

\newcounter{vadik} \newtheorem{theorem}[equation]{Theorem}
\newtheorem{itdefinition}[equation]{Definition}
\newenvironment{definition}{\begin{itdefinition}\rm}{\end{itdefinition}}
\newtheorem{lemma}[equation]{Lemma}
\newtheorem{assume}[equation]{Assumption}  \newtheorem*{cona}{Condition B}
\newtheorem*{conb}{Condition $\mbox{C}_{\kappa}$}

\newtheorem*{theorem*}{Theorem}
\newtheorem{proposition}[equation]{Proposition}
\newtheorem{corollary}[equation]{Corollary}
\newtheorem*{corollary*}{Corollary} 
\newtheorem{remark}[equation]{Remark}
 \newtheorem*{remark*}{Remark}

\newtheorem{itexercise}{Exercise}

\newtheorem*{itexdif8}{*Exercise 8}

\newtheorem*{itexdif9}{*Exercise 9}

\newtheorem*{itexdif10}{*Exercise 10}

\newtheorem*{ithint}{Hint}

\newtheorem*{itinstruction}{Instruction}

\newtheorem{itsolution1}{Solution}

\newtheorem*{itsolution}{Solution}

\date{August 15, 2003; ~Revised April 21, 2004}
\title{Limit theorems for one-dimensional transient random
walks in Markov environments}
\author{
Eddy Mayer-Wolf\thanks{Department of Mathematics, Technion - IIT,
Haifa 32000, Israel (e-mail: emw@tx.technion.ac.il).} \and Alexander
Roitershtein\footnote{Department of Mathematics, Technion - IIT, Haifa
32000, Israel (e-mail: roiterst@tx.technion.ac.il). } \and Ofer
Zeitouni\footnote{Departments of Electrical Engineering and
Mathematics, Technion - IIT, Haifa 32000, Israel, $\mbox{}$ and
Department of Mathematics, University of Minnesota, Minneapolis, MN
55455, USA (e-mail: zeitouni@math.umn.edu)} \footnote{E. M-W.  and
A. R.thank the Department of Mathematics at the University of
Minnesota for its hospitality during a visit in which part of this
work was carried out. The work of O. Z. was partially supported by NSF
grant number DMS-0302230.}}
\begin{document}
\maketitle
\begin{abstract}
We obtain non-Gaussian limit laws for one-dimensional random walk in a
random environment in the case that the environment is a function of a
stationary Markov process. This is an extension of the work of Kesten,
M. Kozlov and Spitzer \cite{kks} for random walks in
i.i.d. environments. The basic assumption is that the underlying
Markov chain is irreducible and either with finite state space or with
transition kernel dominated above and below by a probability measure.
\end{abstract}

\section{Introduction and Statement of Results}
\label{intro}
Let $\Omega=(0,1)^\zz$ and let $\calf$ be the Borel $\sigma-$algebra
on $\Omega.$ A {\em random environment} is an element
$\omega=\{\omega_i\}_{i \in \zz}$ of $\Omega$ distributed according to
a stationary and ergodic probability measure $P$ on $(\Omega,{\mathcal
F})$. The {\em random walk in the environment} $\omega$ is a
time-homogeneous Markov chain $X=$ $\{X_n\}_{n \in {\nn}}$ on $\zz$
governed by the {\em quenched} law
\beq
P_\omega(X_0=0)=1~~~\mbox{and}~~~P_{\omega}(X_{n+1}=j|X_n=i)=\left\{
\begin{array}{ll}
\omega_i&\mbox{if}~j=i+1,\\
1-\omega_i&\mbox{if}~j=i-1.
\end{array}
\right. \feq Let $\left(\zz^\nn,{\mathcal G}\right)$ be the
canonical space for the paths of $\{X_n\},$ i.e. $\mathcal G$ is
the cylinder $\sigma-$algebra. {\em The random walk in random
environment} (RWRE) associated with $P$ is the process
$(X,\omega)$ on the measurable space $\left(\Omega \times \zz^\nn,
{\mathcal F} \otimes {\mathcal G}\right)$ having the {\em
annealed} probability law $\pp=P \otimes P_\omega$ defined by \beq
\pp(F \times G)=\int_F P_\omega(G)P(d\omega),~~~ F \in {\mathcal
F},~ G \in {\mathcal G}. \feq  Since the process learns about the
environment as time passes according to the Bayes rule, $\{X_n\}$
is in general not a Markov chain under the annealed measure $\pp.$
The model goes back to \cite {kozlov73,solomon} and, in physics,
to \cite {chernov67,temkin72}. In this introduction we briefly
discuss some basic results on the one-dimensional RWRE. We refer
the reader to \cite{sznit-topics,notes} for recent comprehensive
surveys of the field.

Recurrence criteria and possible speed regimes for the one-dimensional
RWRE were established by Solomon \cite{solomon} in the case where
$\{\omega_n\}$ is an i.i.d. sequence and carried over to general
ergodic environments by Alili \cite{alili}. Let \beqn
\rho_n&=&\fracd{1-\omega_n}{\omega_n}, \nonumber \\
\label{sf} R(\omega)&=&1+\sum_{n=0}^{+\infty}\rho_0\rho_{-1}
\cdots \rho_{-n}, \feqn $T_0=0,$ and for $n \in \nn,$ \beqn
\label{ti-tau} T_n=\min\{k: X_k \geq
n\}~~~\mbox{and}~~~\tau_n=T_n-T_{n-1}. \feqn $X_n$ is $\as$ transient
if $E_P(\log \rho_0) \neq 0$ and is $\as$ recurrent if $E_P(\log
\rho_0)=0.$ Moreover, if $E_P(\log
\rho_0)<0$ then (see \cite[Sect 2.1]{notes}) $\lim_{n \ra \infty}
\pp(X_n=+\infty)=1,$ $T_n$ are a.s. finite, $\{\tau_n\}$ is a
stationary and ergodic sequence, and we have the following law of
large numbers: \beqn \label{speed} \mbox{v}_P:=\lim_{n \ra
+\infty}\fracd{X_n}{n}=\lim_{n \ra +\infty}\fracd{n}{
T_n}=\fracd{1}{\ee(\tau_1)}= \fracd{1}{2E_P(R)-1},~~~\pp-\as \feqn
Thus, the transient walk $X_n$ has a deterministic speed
$\mbox{v}_P=\lim_{n \ra \infty}X_n/n$ which may be zero.

Solomon's law of large numbers for the transient walks in i.i.d.
environment was completed by limit laws in the work of  Kesten, M.
Kozlov, and Spitzer \cite{kks}. The limit laws for the RWRE $X_n$
are deduced in \cite{kks} from stable limit laws for the hitting
times $T_n,$ and the index $\kappa$ of the stable distribution is
determined by the condition \beq E_P(\rho_0^\kappa)=1.\feq In
particular, under certain conditions the central limit theorem
holds with the standard normalization $\sqrt{n},$ and this case
was extended to stationary and ergodic environments by Alili
\cite{alili}, Molchanov \cite{molchanov} and Zeitouni \cite[Sect
2.2]{notes}, see also Bremont \cite{bremont}.

In this paper we obtain limit laws for $X_n$ for environments which
are pointwise transformations of a stationary ergodic Markov process
which satisfies Assumption \ref{measure1} below.  These laws are
related to stable laws of index $\kappa \in (0,2],$ where, under the
assumptions below, $\kappa$ is determined by
\beqn \label{kappa239} \Lambda(\kappa)=0\,,\quad\mbox{\rm
where}\,\, \Lambda(\beta):=\lim_{n \ra \infty} \fracd{1}{n} \log E_P
\left(\proda_{i=0}^{n-1} \rho_i^\beta\right). \feqn More precisely:
\\{\bf Basic setup:} On a state space $\cals$ equipped with a
countably generated $\sigma-$algebra $\calt,$ let $\{x_n\}_{n
\in \zz}$ be a stationary Markov chain, such that
$\omega_{-n}=\omega(x_n)$ (and hence $\rho_{-n}=\rho(x_n)$) for
measurable functions $\rho,\omega:\cals \ra \rr.$ We denote by
$H(x,\cdot)$ the transition probability measure of $(x_n),$ by
$\pi$ its stationary probability measure, and use the notation
$H(x,y)$ to denote $H(x,\{y\})$ for a single state $y \in \cals.$
With $P_x$ denoting the law of the Markov chain with $x_0=x$, the
reader should not confuse $P_x$ and $P_\omega$.

We shall say that the process $\log \rho_{-n}$ is $\alpha$-arithmetic
(c.f. [22,2]) if $\alpha >0$ is the largest number for which there
exists a measurable function $\gamma :\cals \ra [0,\alpha)$\ such
that 
\beq
P(\log\rho_0\in\gamma(x_{-1})\!-\!\gamma(x_0)\!+\!\alpha\mathbb{Z})=1,
\ \ \ P-{\rm a.s.} 
\feq
The process will be said to be non-arithmetic if no such $\alpha$
exists.

\begin{assume} \label{measure1}
\item [(A1)] Either \beqn \label{finite} \cals~\mbox{\em is a
finite set and the Markov chain}~(x_n)~\mbox{\em is irreducible},
\feqn or, there exist a constant $c_r \geq 1$ and a probability
measure $\psi$ on $(\cals,\calt)$ such that for some $m \in \nn,$
\beqn \label{lem-lem} c^{-1}_r\psi(A) < H^m(x,A)<c_r \psi(A)~~~
\forall x \in \cals,A \in \calt, \feqn where the kernel $H^n(x,A)$
is defined inductively by $H^0(x,A)=\odin_{A}(x)$ for all $x \in \cals,
A\in \calt$ and $H^n(x,A)=\int_\cals H^{n-1}(x,dy)H(y,A),$ $n \geq
1.$ \item [(A2)] $P(\epsilon<\omega_0<1-\epsilon)=1$ for some
$\epsilon \in (0,1/2).$ 
\item [(A3)] $\limsup_{n \ra \infty} \fracd{1}{n} \log E_P
\left(\prod_{i=0}^{n-1} \rho_i^\beta \right)<0$ and $\limsup_{n \ra 
\infty} \fracd{1}{n} \log E_P \left( \prod_{i=0}^{n-1} 
\rho_i^{\beta'}\right) \geq 0$ for some constants $\beta>0$ and 
$\beta'>0.$
\item [(A4)] $\log \rho_0$ is non-arithmetic in the sense defined above.
\end{assume}
Note that condition {\em (A1)} refers to the underlying Markov chain
$(x_n),$ whereas conditions {\em (A2)}--{\em (A4)} refer to $\omega$
itself. Assumption \eqref{finite} is not a particular case of
assumption \eqref{lem-lem} since under \eqref{finite} the Markov
chain $(x_n)$ may be periodic. Under {\em (A1)}, the environment
$\omega$ is an ergodic sequence (see e.g. \cite[p.  338]{durrett} or
\cite[Theorem 6.15]{nummelin}). Condition {\em (A3)} guarantees, by
convexity, the existence of a unique $\kappa$ in
\eqref{kappa239}. Indeed it will be shown later that the $\limsup$ is
in fact a $\lim.$ It also follows from {\em (A3)}, by Jensen's
inequality, that $E_P(\log
\rho_0)<0,$ so that $X_n$ is transient to the right. For future
reference we denote \beqn
\label{crho} c_\rho=\fracd{1-\epsilon}{\epsilon},
\feqn and note that by the ellipticity condition {\em (A2)},
$P(c_\rho^{-1}<\rho_0<c_\rho)=1.$

For $\kappa \in (0,2]$ and $b >0$ we denote by $\call_{\kappa,b}$
the stable law of index $\kappa$ with the characteristic function
\beqn \label{kappa-law} \log \widehat
\call_{\kappa,b}(t)=-b|t|^\kappa\left(1+i\fracd{t}{|t|}
f_\kappa(t)\right), \feqn where $f_\kappa(t)=-\tan
\fracd{\pi}{2}\kappa$ if $\kappa \neq 1,$ $f_1(t)=2/\pi \log t.$
With a slight abuse of notation we use the same symbol for the
distribution function of this law. If $\kappa<1,$
$\call_{\kappa,b}$ is supported on the positive reals, and if
$\kappa \in (1,2],$ it has zero mean \cite[Chapter
1]{samor-taqqu}. Our main result is:
\begin{theorem}
\label{main-markov} Let Assumption \ref{measure1} hold. Then there
is a unique $\kappa>0$ such that \eqref{kappa239} and the
following hold for some $b>0:$ \item[(i)] If $\kappa \in (0,1),$
then $\lim_{n \ra \infty}\pp\left(n^{-\kappa}X_n \leq \mathfrak{z}
\right)=1-\call_{\kappa,b}(\mathfrak{z}^{-1/\kappa}),$ \item[(ii)]
If $\kappa=1,$ then $\lim_{n \ra \infty}\pp\bigl(n^{-1}(\log
n)^2(X_n-\delta(n)) \leq \mathfrak{z}
\bigr)=1-\call_{1,b}(-\mathfrak{z}),$ for suitable $A_1>0$ and
$\delta(n) \sim (A_1\log n)^{-1} n,$ \item[(iii)] If $\kappa \in
(1,2),$ then $\lim_{n \ra \infty}
\pp\left(n^{-1/\kappa}\left(X_n-n \mbox{\em v}_P\right) \leq
\mathfrak{z}\right)=1- \call_{\kappa,b}(-\mathfrak{z}).$
\item[(iv)] If $\kappa=2,$ then $\lim_{n \ra \infty} \pp\left((n
\log n)^{-1/2} (X_n-n \mbox{\em v}_P)\leq
\mathfrak{z}\right)=\call_{2,b}(\mathfrak{z}).$
\end{theorem}
In the setup of Theorem \ref{main-markov} it is not hard to check, and
follows e.g. from \cite[Theorem 2.2.1]{notes}, that the standard CLT
holds if $\kappa>2.$

As in \cite{kks}, stable laws for $X_n$ follow from stable laws
for the hitting times $T_n$, and we direct our efforts to
obtaining limit laws for the latter. We have:
\begin{proposition}
\label{markov-tau} Let Assumption \ref{measure1} hold. Then there
is a unique $\kappa>0$ such that \eqref{kappa239} and the
following hold for some $\tilde b>0:$ \item[(i)] If $\kappa \in
(0,1),$ then $\lim_{n \ra \infty}\pp\left(n^{-1/\kappa}T_n \leq t
\right)=\call_{\kappa,\tilde b}(t),$ \item[(ii)] If $\kappa=1,$
then $\lim_{n \ra \infty}\pp\bigl(n^{-1}(T_n-nD(n)) \leq t
\bigr)=\call_{1,\tilde b}(t),$ for suitable $c_0>0$ and $D(n) \sim
c_0 \log n,$\item[(iii)] If $\kappa \in (1,2),$ then $\lim_{n \ra
\infty} \pp\left(n^{-1/\kappa}\left(T_n-n \mbox{\em
v}_P^{-1}\right) \leq t\right)=\call_{\kappa,\tilde b}(t).$
\item[(iv)] If $\kappa=2,$ then $\lim_{n \ra \infty} \pp\left((n
\log n)^{-1/2}(T_n-n \mbox{\em v}_P^{-1})\leq
t\right)=\call_{2,\tilde b}(t).$
\end{proposition}
The proof that Theorem \ref{main-markov} follows from Proposition
\ref{markov-tau} is the same as in the i.i.d. case, and is based
on the observation that for any positive integers $\eta,\zeta,n$
\beqn \label{unif} \{T_\zeta \geq n\} \subset \{ X_n \leq \zeta\}
\subset \{T_{\zeta+\eta} \geq n\} \bigcup \{\inf_{k \geq
T_{\zeta+\eta} } X_k -(\zeta+\eta) \leq -\eta\}. \feqn Because the
random variables $\inf_{k \geq T_{\zeta+\eta} } X_k -(\zeta+\eta)$
and $\inf_{k \geq 0} X_k$ have the same annealed distribution, the
probability of the last event in \eqref{unif} can be made
arbitrary small uniformly in $n$ and $\zeta$ by fixing $\eta$
large (since the RWRE $X_n$ is transient to the right). For
$\kappa=1,$ the rest of the argument is detailed in \cite[pp.
167--168]{kks}, where no use of the i.i.d. assumption for $\omega$
is made at that stage, and a similar argument works for all
$\kappa \in (0,2].$ All of our work in the sequel is directed
toward the proof of Proposition \ref{markov-tau}.

Following \cite{kks}, the analysis of $T_n$ is best understood in
terms of certain regeneration times $\nu_n$, with excursion counts
between regenerations forming a branching process $Z_n$ with
immigration in a random environment (see Section \ref{branch} for
precise definitions). In the i.i.d. setup, the total population of
the branching process between regenerations, denoted $W_n$, forms
an i.i.d. sequence, and much of the work in \cite{kks} is to
establish accurate enough tail estimates on them to allow for the
application of the i.i.d. stable limit law for partial sums of
$W_n.$ The limit laws for $T_n$ then easily follow from those for
$W_n.$

In our case, the sequence $W_n$ a-priori is not even stationary.
However, using the regeneration property of the underlying Markov
chain $(x_n)$ (see Section \ref{proper}), we introduce in Section
\ref{branch} modified regeneration times $\bar \nu_n$ (a random
subsequence of $\nu_n$) such that the total population of the
branching process between times $\bar \nu_n$ and $\bar \nu_{n+1},$
denoted by $\ol W_{n+1},$ is a one-dependent stationary sequence.
This sequence is i.i.d. if either \eqref{lem-lem} with $m=1$ or
\eqref{finite} hold. Again following the proof in \cite{kks}, we
obtain tails estimates for the random variables $\ol W_{n+1}$
yielding the stable limit laws for $T_n$ stated in Proposition
\ref{markov-tau}. Similarly to the i.i.d. case, the key to the
proof is the derivation of tails estimates obtained in Section
\ref{proof-ar} for the random variable $R$ defined in \eqref{sf}.

We conclude the introduction with a characterization of the speed
$\mbox{v}_P$ under Assumption \ref{measure1}, which will not be
used in the sequel. Recall that $\rho_{_n}=\rho(x_n)$ for a
measurable function $\rho: \cals \ra \rr.$ If $\kappa \leq 1,$
then $\mbox{v}_P=0,$ and if $\kappa>1,$ then $\mbox{
v}_P^{-1}=E_P\bigl(\rho(x_0) \xi(x_0)\bigr),$ where the function
$\xi: \cals \ra (0,\infty)$ is the unique positive and bounded
solution of the equation \beqn \label{markov-speed}
\xi(x)=\int_\cals H(x,dy)\rho(y)\xi(y)+1+1/\rho(x). \feqn This
formula is essentially due to Takacs \cite{takacs}, who considered
finite-state Markov environments. The proof in the general Markov case
is included at the end of Section \ref{proper}.

The rest of the paper is organized as follows. Section
\ref{proofs}, divided into three subsections, contains the proof
of Theorem \ref{main-markov}, except for the proofs of two
propositions which are deferred to the Appendix. In Subsection
\ref{proper} some basic properties of Markov chains that satisfy
Assumption \ref{measure1} are described. In particular, Condition
B is introduced and shown to hold under Assumption \ref{measure1}.
In Subsection \ref{branch}, Condition $\mbox{C}_\kappa$ is
introduced and Proposition \ref{markov-tau} is derived from it and
Condition B, making use of the above mentioned branching process
and a regeneration structure it possesses. Finally, Subsection
\ref{proof-ar} is devoted to the proof that Condition
$\mbox{C}_\kappa$ holds under Assumption \ref{measure1}.
\section{Proofs}
\label{proofs}
\subsection{Some properties of the underlying Markov chain and their 
consequences}
\label{proper} We summarize here, using the framework of the
Athreya-Ney and Nummelin theory of positive recurrent kernels (cf.
\cite{athreya-ney,athreya-ney1,nummelin}), some properties of the
Markov chain $(x_n)$ that follow from Assumption \ref{measure1}.
The main objectives here are to introduce the regeneration times
$N_k$ and to obtain the Perron-Frobenius type Lemmas \ref{urg} and
\ref{urg1}. One immediate consequence of these lemmas is that
Condition B introduced subsequently is satisfied under Assumption
\ref{measure1}.

First, we define a sequence of regeneration times for the Markov
chain $(x_n).$ If \eqref{finite} holds, let $x^* \in \cals$ be any
(recurrent) state of the Markov chain $(x_n)$ and pick any $r \in
(0,1).$ Let $(y_n)_{n \in \zz}$ be a sequence of i.i.d.  variables
independent of $(x_n)$ (in an enlarged probability space if
needed) such that $P(y_0=1)=r$ and $P(y_0=0)=1-r,$ and let \beq
N_0=0,~~~N_{n+1}=\min\{k>N_n: x_n=x^*,~y_n=1\},~n \geq 0. \feq
Then, the blocks $\bigl(x_{_{N_{_n}}},x_{_{N_{_n}+1}}, \ldots,
x_{_{N_{_{n+1}}-1}}\bigr)$ are independent, and $x_{_{N_{_n}}}$
are identically distributed for $n \geq 1.$ Note that between two
successive regeneration times, the chain evolves according to the
sub-stochastic Markov kernel $\Theta$ defined by \beqn
\label{kernels1} H(x,y)=\Theta(x,y)+r \odin_{\{y=x^*\}}
H(x,y),\feqn that is \beqn \label{between}
P_x(x_1=y,N_1>1)=\Theta(x,y).\feqn If \eqref{lem-lem} holds, then
the random variables $N_k$ can be defined by the following
procedure (see \cite{athreya-ney, nummelin} and \cite{asmus}).
Given an initial state $x_0,$ generate $x_m$ as follows: with
probability $r<c_r^{-1}$ distribute $x_m$ over $\cals$ according
to $\psi$ and with probability $1-r$ according to $1/(1-r)\cdot
\Theta(x_0,\cdot),$ where the kernel $\Theta(x,\cdot)$ is defined
by \beqn \label{kernels} H^m(x,A)=\Theta(x,A)+r\psi(A),~~~x \in
\cals,A \in \calt. \feqn Then, (unless $m=1$) sample the segment
$\bigl(x_1,x_2,\ldots, x_{m-1}\bigr)$ according to the chain's
conditional distribution, given $x_0$ and $x_m.$ Generate $x_{2m}$
and $x_{m+1},x_{m+2},\ldots,x_{2m-1}$ in a similar way, and so on.
Since the ``$r$-coin" is tossed each time independently, the event
``the next move of the chain $(x_{mn})_{n \geq 0}$ is according to
$\psi$" occurs i.o. Let $N_0=0$ and $\{N_k\}_{k \geq 1}$ be the
successful times of its occurrence multiplied by $m.$ By
construction, the blocks $\bigl(x_{_{N_{_n}}},x_{_{N_{_n}+1}},
\ldots, x_{_{N_{_{n+1}}-1}}\bigr)$ are one-dependent (if $m=1$
they are actually independent), and for $n \geq 1$ they are
identically distributed ($x_{_{N_{_n}}}$ is distributed according
to $\psi$).

Let us summarize the most important property of the regeneration
times $N_n$ as follows. For $n \geq 0,$ let \beqn
\label{di-blocks} D_n=\bigl(x_{_{N_{_n}}},x_{_{N_{_n}+1}}, \ldots,
x_{_{N_{_{n+1}}-1}}\bigr).\feqn Then:
\begin{itemize}
\item The random blocks $D_n$ are identically distributed for $n
\geq 1.$ \item If \eqref{finite} or \eqref{lem-lem} with $m=1$
hold, $D_n$ are independent for $n \geq 0.$ \item If
\eqref{lem-lem} holds with $m>1,$ $D_n$ are one-dependent for $n
\geq 0.$
\end{itemize}
In both cases under consideration (either of \eqref{finite} or of
\eqref{lem-lem}), there exist constants $l,\delta>0,$ such that
(cf. \cite{athreya-ney}) \beqn \label{as} \inf_{x \in
\cals}P_x(N_1 \leq l)>\delta>0. \feqn The regeneration times $N_n$
will be used in Section \ref{branch} for the construction of an
auxiliary sequence $\ol W_n$ of stationary and one-dependent
random variables playing a central role in the proof of
Proposition \ref{markov-tau}.

We now turn to a Perron-Frobenius type theorem for positive finite
kernels, having in mind applications to the kernels of the form
$K(x,A)=E_x\left(\prod_{i=0}^n \rho_{-i}^\beta; x_n \in A\right).$
In the following two lemmas, we consider separately the cases of
non-finite (assumption \eqref{lem-lem}) and finite (assumption
\eqref{finite}) state space $\cals.$ In particular, the properties
of the positive kernels described in these lemmas imply Condition
B introduced below and are essential for the proof of the crucial
Proposition \ref{sb}.

Let $B_b$ be the Banach space of bounded measurable real-valued
functions on $(\cals,\calt)$ with the norm $\|f\|=\sup_{x \in
\cals}|f(x)|.$ A positive and finite kernel $K(x,A)$ (a measurable
function of $x$ for all $A \in \calt$ and a finite positive
measure on $\calt$ for all $x \in \cals$) defines a bounded linear
operator on $B_b$ by setting $Kf(x)=\int_\cals K(x,dy)f(y).$ We
denote by $r_{_K}$ the spectral radius of the operator
corresponding to the kernel $K,$ that is $$r_{_K}=\lim_{n \ra \infty}
\sqrt[n]{\|K^n \odin\|}=\lim_{n\ra \infty}
\sqrt[n]{\|K^n\|_{B_b\to B_b}},$$ where $\odin(x) \equiv 1.$

Although the results stated in the following lemma are certainly
well-known and appear elsewhere, their proofs are provided for the
sake of completeness.
\begin{lemma}
\label{urg} Let $K(x,A)$ be a positive kernel on $(\cals,\calt)$
such that for some constant $c \geq 1$ and probability measure
$\psi,$ \beqn \label{ka-kernel} c^{-1} \psi(A) \leq K(x,A) \leq c
\psi(A),~~~\forall x \in \cals,~A \in \calt. \feqn  Then,
\item[(a)] There exists a function $f\in
B_b$ such that $\inf_x f(x)>0$ and $Kf=r_{_K} f.$ There exists a constant
$c_K \geq 1$ such that $c_K^{-1} r_{_K}^n \leq K^n \odin
\leq c_K r_{_K}^n$ for all $n \in \nn.$ \item[(b)]
If $K=K^m_1$ for a positive finite kernel $K_1(x,A)$ and some $m \in
\nn,$ then $r_{_{K_1}}=r_{_K}^{1/m}$ and there exists a function $f_1\in
B_b$ such that $\inf_x f_1(x)>0$ and $K_1f_1=r_{_{K_1}}^{1/m} f_1.$
\end{lemma}
\begin{proof}
\item[{\em (a)}] The existence of a function $f: \cals \ra
(0,\infty)$ and a constant $\lambda >0$ such that $K f =\lambda f$
follows from the Example in \cite[p. 96]{nummelin}. It follows
from \eqref{ka-kernel} that $f(x)$ is bounded away from zero and
infinity, i.e. $c_{_K}^{-1}\leq f(x)\leq c_{_K}$ for some
$c_{_K}>0.$ Hence, for any $n>0,$ $K^n \odin< c_{_K} K^n f =c_{_K}
\lambda^n f <c_{_K} ^2 \lambda ^n .$ Similarly, $K^n
\odin>c_{_K}^{-2} \lambda^n .$ That is, $\lambda =r_{_K}.$
\item[{\em (b)}] Set $f_1=\sum_{j=0}^{m-1} (1/r_{_K})^{j/m}
K^j_1f.$
\end{proof}
The finite-state counterpart of the previous lemma is stated as
follows:
\begin{lemma}
\label{urg1} Let $\cals=\{1,2,\ldots,n\}$ and $K(i,j)$ be an
irreducible $n \times n$ matrix with nonnegative entries. For some
constants $r \in (0,1)$ and $j^* \in \{1,\ldots,n\}$ define the matrix
$\witi \Theta(i,j)$ by \beqn \label{theta-kernel-1} K(i,j)=\witi
\Theta(i,j)+r\odin_{\{j=j^*\}}K(i,j),~~~ 1 \leq i,j
\leq n. \feqn
Then,
\item [(a)]  Assertion (a) of Lemma \ref{urg} holds for the matrix $K.$
\item[(b)] There exists a function $g \in B_b$ such that
$\inf_x g(x)>0$ and $\witi \Theta g=r_{\witi \Theta} g.$ \item[(c)]
$r_{_{\witi \Theta}}\in (0,r_{_K}).$
\end{lemma}
\begin{proof}
Since $\witi \Theta$ and $K$ have the same adjacency matrices
($K(i,j)=0$ iff $\witi \Theta(i,j)=0$), $\witi \Theta$ is
irreducible as well. Assertions of {\em (a)} and {\em (b)} follow
then from the Perron-Frobenius theorem.  Clearly $r_{_{\witi \Theta}} \leq
r_{_K}.$ Since $r_{_K} f \geq \witi \Theta f,$ the equality
$r_{_{\witi \Theta}}=r_{_K}$ would imply \cite[Theorem
5.1]{nummelin} that $f=g$ and $\witi \Theta f =r_{_K} f =K f ,$ that
is impossible since $f>0$ everywhere. Hence $r_{_{\witi \Theta}}<
r_{_K}.$
\end{proof}
Since for any $\beta \geq 0,$ \beqn \label{hbeta}E_x
\left(\proda_{k=0}^{n-1}~(\rho_{-k})^\beta\right)= \rho(x)^\beta
H_\beta^{n-1} \odin(x),\feqn where
$H_\beta(x,dy)=H(x,dy)\rho(y)^\beta,$ it follows from Lemmas
\ref{urg} and \ref{urg1} that for some constant $c_\beta \geq 1$
which depends on $\beta$ only, \beqn \label{h-beta} c_\beta^{-1}
r_\beta^n \leq E_x\left(\proda_{k=0}^{n-1}~(\rho_{-k})^\beta\right)
\leq c_\beta r_\beta^n,~~~x \in \cals,~n \in \nn, \feqn where
$r_\beta=r_{_{H_\beta}}.$ Therefore, the following Condition B is
satisfied under Assumption \ref{measure1}. With future applications in
mind, we make the formulation suitable for non-Markovian ergodic
environments. Let \beqn \label{f-zero}
\calf_0=\sigma(\omega_n: n >0)\feqn be the $\sigma-$algebra
generated by the ``past'' of the sequence $\{\omega_{-n}\}.$
\begin{cona}
\label{a239} $\{\omega_{-n}\}$ is a stationary and ergodic
sequence such that \item[(B1)] Ellipticity condition:
$P(\epsilon<\omega_0<1-\epsilon)=1$ for some $\epsilon \in
(0,1/2).$ \item [(B2)] For any $\beta >0,$  \beqn
\label{conb-cite} \lim_{n \ra \infty} \fracd{1}{n} \log E_P
\left(\proda_{k=0}^{n-1}~\rho_{-k}^\beta ~\bigl|
\calf_0\right)=\Lambda(\beta),\; \as, \feqn with uniform (in
$\omega$) rate of convergence, with $\Lambda(\beta)$ as in
\eqref{kappa239}. Further, there exists a unique $\kappa>0$ such
that $\Lambda(\kappa)=0,$ and $\Lambda(\beta)(\beta-\kappa) \geq
0$ for all $\beta>0.$
\end{cona}
The last statement follows since $\Lambda(\beta)$ is a convex
function of $\beta$ in $[0,\infty),$ taking both negative and
positive values  by Assumption {\em (A3)}, with $\Lambda(0)=0.$
\\
\centerline{} \\ We conclude this subsection with the proof of
\eqref{markov-speed}. It follows from \eqref{speed}, \eqref{sf}
and \eqref{h-beta} that $\mbox{v}_P=0$ for $\kappa \leq 1.$ Assume
that $\kappa>1$ and consider the following decomposition for the
hitting time $\tau_1$ defined in \eqref{ti-tau}): \beq \tau_1={\bf
1}_{\{X_1=1\}}+{\bf 1}_{\{X_1=-1\}}(1+\tau_0''+\tau_1'), \feq
where $1+\tau_0''$ is the first hitting time of 0 after time 1,
and $1+\tau_0''+\tau_1'$ is the first hitting time of 1 after time
$1+\tau_0''$.  Taking expectations in both sides of the equation
(first for a fixed environment and then integrating over the set
of environments) gives \beq
\ee(\tau_1|x_0=x)=1+\rho(x)\left(1+\ee(\tau_0''|x_0=x) \right).
\feq Since $\ee(\tau_0''|x_0=x)=\ee(\tau_1|x_1=x) =\int_\cals
\ee(\tau_1|x_0=y)H(x,dy),$ we obtain that the function
$\xi(x):=\ee(\tau_1|x_1=x)/\rho(x)$ solves equation
\eqref{markov-speed}. Recalling the operator $H_1:f(x) \ra$ $ \int_\cals
H(x,dy)\rho(y)f(y)$ acting on $B_b,$ it follows from identity
\eqref{h-beta} and Condition B, that its spectral radius is strictly
less than one, and a simple truncation argument (by
\eqref{markov-speed}, $\xi_M \leq H_1 \xi_M+1+1/\rho,$ where
$\xi_M(x):=\ee\left(\min\{\tau_1,M\}|x_1=x\right)/\rho(x)$ for a
constant $M>0$) shows that $\xi(x)$ is a bounded function of $x,$
yielding that $\ee(\tau_1)=E_P\bigl(\rho(x_0)\xi(x_0)\bigr).$ This
implies \eqref{markov-speed} by \eqref{speed} (Lemmas 2.1.11 and
2.1.17 in \cite{notes}).
\subsection{The branching model and its regeneration structure}
\label{branch} We consider here a branching process $\{Z_n\}$ in
random environment with immigration closely related to the RWRE (see
e.g., \cite {alili,kks, notes}). The random variables $T_n$ are
associated by \eqref{ti-zed} to the partial sums of the branching
process $Z_n.$ This leads us naturally to the variables $\ol W_n,$
defined in \eqref{wol-def}, which are random partial sums of $Z_n.$
The aim in introducing the branching process is to transform the limit
problem of $T_n$ into a limit problem for the partial sums of the
sequence $\ol W_n,$ which turns out to be a stationary and
one-dependent sequence in a stable domain of attraction.

Let \beq U_i^n=\# \{k<T_n:X_k=i,~X_{k+1}=i-1\},~~~i,n \in \zz,
\feq the number of moves to the left from site $i$ up to time
$T_n.$ Then \beqn \label{ti-zed} T_n=n+2\sum\limits_{i=-\infty}^n
U_i^n. \feqn When $U^n_n=0,U^n_{n-1}, \ldots, U^n_{n-i+1}$ and
$\omega_n, \omega_{n-1} \ldots,\omega_{n-i}$ are given,
$U^n_{n-i}$ is the sum of $U^n_{n-i+1}+1$ i.i.d. geometric random
variables that take the value $k$ with probability
$\omega_{n-i}(1-\omega_{n-i})^k,$ $k=0,1,\ldots$ Assuming that the
RWRE is transient to the right we have: \beqn \label {spent-neg}
\sum\limits_{i \leq 0} U_i^n \leq ~\mbox{total time spent by}
~\{X_t\} ~\mbox{in}~ (-\infty;0]<\infty~\as\feqn Therefore, in
order to prove the limit laws for $T_n$ it is sufficient to prove
the corresponding result for the sums $\sum_{i=1}^n U_i^n.$ These
sums have the same distribution as \beqn \label{zed-sums}
\sum\limits_{k=0}^{n-1}Z_k, \feqn where $Z_0=0,Z_1,Z_2,\ldots$
forms a branching process in random environment with one immigrant
at each unit of time.

Without loss of generality, we shall extend the underlying sample
space $\left(\Omega \times \zz^\nn\right)$ to $\left(\Omega \times
\Upsilon\right),$ where $\Upsilon$ is large enough to fit not only
the random walk but also the branching process, and assume that
$P_\omega$ (and hence $\pp$) is suitably extended.

Thus, when $\omega$ and $Z_0, \ldots,Z_n$ are given, $Z_{n+1}$ is
the sum of $Z_n+1$ independent variables $V_{n,0},V_{n,1},\ldots,
V_{n,Z_n}$ each having the geometric distribution \beqn
\label{bigv} P_\omega
\{V_{n,j}=k\}=\omega_{-n}(1-\omega_{-n})^k,~k=0,1,2,\ldots \feqn
Extending \eqref{f-zero}, let for $n \in \nn,$ \beqn
\label{sigmag}
\calf_n=\sigma\left(Z_0,Z_1,Z_2,\ldots,Z_{n-1},Z_n\right) \vee
\sigma\left(\omega_j: j>-n \right),\feqn that is, the
$\sigma$-algebra generated by the branching process
$\{Z_i\}_{i=0}^n$ and the environment
$\{\omega_i\}_{i=-n+1}^\infty$ before time $n.$

As in \cite{kks}, the random variables \beq
\nu_0=0,~~~\nu_n=\min\{k>\nu_{n-1}:Z_k=0\} \feq are the successive
stopping times at which the population becomes extinct, and the
variables \beq W_n=\suml_{k=\nu_{n-1}}^{\nu_n-1}Z_k \feq measure
the total number of individuals born between two such extinction
times.

Recall the definition of the $\sigma-$algebra $\calf_0$ given in
\eqref{f-zero}. The proof of the following proposition, which is a modification
of Lemma 2 in \cite{kks} adapted to non-i.i.d. environments,
is included in Appendix \ref{moment-proof}.
\begin{proposition}
\label{moment} Assume that Condition B holds. Then, there exist
$C_1,~C_2>0$ such that $P-\as,$ $\pp( \nu_1>n |\calf_0) \leq
C_1e^{-C_2n},$ for any $n>0.$
\end{proposition}
The following corollary is immediate since $C_1,C_2$ above are
deterministic.
\begin{corollary}
\label{momentc} Assume that Condition B holds. Then, with
probability one, $\pp( \nu_{j+1}-\nu_j>n |\calf_{\nu_j}) \leq
C_1e^{-C_2n},$ for any $j \geq 0$ and $n>0,$ where the constants
$C_1,~C_2>0$ are the same as in Proposition \ref{moment}.
\end{corollary}
Let $\{N_k\}_{k=0}^\infty$ be the sequence of successive
regeneration times for the chain $(x_n)$ defined in Section
\ref{proper}, let $\bar \nu_0=0,$ and for $n \geq 0$ define the
stopping times: \beqn \label{nubar} \bar \nu_{n+1}=\inf\{k> \bar
\nu_n: k=\nu_i=N_j~\mbox{for some}~i,j>0\},\feqn and the random
variables \beqn \label{wol-def} \ol W_{n+1}=\sum_{k=\bar
\nu_n}^{\bar \nu_{n+1}-1} Z_k. \feqn By construction of the random
times $N_n,$ the segments of the environment between $\bar \nu_n$
and $\bar \nu_{n+1}-1$ are one-dependent (see \eqref{di-blocks}
and the subsequent summary), and hence the variables $\{\ol
W_n\}_{n \geq 1}$ form a one-dependent sequence, which is even
independent if either \eqref{finite} or \eqref{lem-lem} with $m=1$
hold.
\begin{lemma} \label{st17} Let Assumption \ref{measure1} hold.
Then, \item[(a)] The distribution of $\bar\nu_1$, conditioned on
the ``past'' has exponential tails: there exist $K_1,$ $K_2>0$
such that $P-\as,$ \beqn \label{moment3} \pp(
\bar{\nu}_1>n|\calf_0 ) \leq K_1e^{-K_2 n},~~~\forall n>0, \feqn
and, more generally, \beqn \label{3moment}
\pp( \bar \nu_{j+1}-\bar
\nu_j>n|\calf_{\bar \nu_j} ) \leq K_1e^{-K_2 n}\feqn
for any $j \geq 0.$
\item[(b)] The law of large numbers holds for $\bar \nu_n:$ $\pp
\left(\lim_{n \ra \infty} \fracd{\bar \nu_n}{n}=\mu\right)=1,$ where
$\mu=\ee(\bar \nu_2-\bar \nu_1)>0.$ \item[(c)] The central limit
theorem holds for $\bar \nu_n:$ there exists a constant $b
>0$ such that the law of $(\bar \nu_n-n\mu)/\sqrt{n}$ converges to
$\call_{2,b}.$
\end{lemma}
\begin{proof}
\item[{\em (a)}]   Clearly, it is sufficient to prove
\eqref{moment3}, since the constants $K_1$ and $K_2$ are
deterministic. Let $F_1=\{Z_1=0\},$ and for $2 \leq j \leq l,$
where $l$ is defined in \eqref{as}, \beq F_j =
\{Z_1=Z_2=\ldots=Z_{j-1}=1,Z_j=0\}, \feq and \beq \cals_j=\{x \in
\cals: P_x(N_1=j)>\delta/l\}. \feq Then $\bigcup_{j=1}^l
\cals_j=\cals,$ and we have for $x \in \cals_j$ : \beq
\pp\bigl(\nu_1=N_1 \leq l|x_0=x) &\geq& \pp \bigl(F_j \cap
\{N_1=j\} |x_0=x \bigr) =\\&=&  P_x(N_1=j)\pp \bigl(F_j |x_0=x,
N_1=j \bigr) \geq \fracd{\delta}{l} \pp\bigl(F_j |x_0=x, N_1=j
\bigr). \feq Using the ellipticity condition {\em (A2)}, we obtain
that $P-\as,$ $P_\omega (F_1)=\omega_0 \geq \epsilon,$ and for $2
\le j \leq l,$ \beq P_\omega (F_j )=\omega_0(1-\omega_0)
\prod_{k=1}^{j-2}
\bigl(2\omega_{-k}^2(1-\omega_{-k})\bigr)\omega_{-j+1}^2 \geq
2^{j-2}\epsilon^{2j-1}(1-\epsilon)^{j-1} \geq
\epsilon^{2l}(1-\epsilon)^{l-1}, \feq implying that
$\pp\bigl(\nu_1=N_1 \leq l|x_0=x) \geq \delta/l \cdot
\epsilon^{2l}(1-\epsilon)^{l-1}>0$ for $P-$almost every $x \in
\cals.$ Thus, in view of Corollary \ref{momentc}, $\bar \nu_1$ is
stochastically dominated by a sum of a geometric random number of
i.i.d. variables with exponential tails, yielding \eqref{moment3}.
We note in passing that, in view of the uniform bounds in the
proof above, the same argument yields uniform exponential tails
for the distribution of $\bar \nu_{i+1}-\bar\nu_i$ conditioned on
$\sigma\{\omega_j, j> -\bar\nu_i\}$. \item[{\em (b)}] Follows from
\eqref{moment3} and the ergodic theorem, since
$\bar\nu_{n+1}-\bar\nu_{n},$ $n \geq 1,$ are one-dependent
identically distributed variables. \item[{\em (c)}] Follows e.g.
from the CLT for stationary and uniformly mixing sequences
\cite[p. 427]{durrett}.
\end{proof}
Recall the function $R(\omega)$ defined in \eqref{sf}. We shall prove
in Subsection \ref{proof-ar} that under Assumption \ref{measure1} the
following condition holds for some $\kappa>0.$
\begin{conb}
\label{conba} There exists a strictly positive random variable
function $K(\omega)$ such that for some positive constants
$K_3,K_4,t_c$ the following hold $P-\as:$ \beqn \label{bound}
t^\kappa P\left(R>t | \calf_0\right)>K_3~~\forall
t>t_c~~~\mbox{and}~~~t^\kappa P\left(R>t |
\calf_0\right)<K_4~~\forall t>0, \feqn \beqn \label{gvul} \lim_{t
\ra \infty}t^\kappa P\left(R>t |\calf_0\right)=K(\omega). \feqn
\end{conb}
It follows from \eqref{bound} and \eqref{speed} that the case $\kappa
\leq 1 $ corresponds to zero speed, and the case $\kappa>1$ to a
positive speed. Note that if Condition $\mbox{C}_{\bar \kappa}$ and
Condition B hold simultaneously, then $\bar \kappa=\kappa.$

For $n \geq 1$ let \beq \witi W_n=\sum_{j=1}^n \ol W_j,\feq where
the random variables $\ol W_j$ are defined in \eqref{wol-def}. The
next proposition is an analogue of \cite[Lemma 6]{kks} for
non-i.i.d environments and is applicable for non-Markov
environments too. 
\begin{proposition}
\label{st39} Assume  Conditions B and $\mbox{C}_{\kappa}$. Then,
for any $n \geq 1$ there exist constants $t_n,L_n,J_n>0$ and a
strictly positive random variable $\witi K_n(\omega)$ such that
the following hold $P-\as:$ \beqn \label{plus} J_n<t^\kappa
\pp\left(\witi W_n>t|\calf_0\right),~\forall
t>t_n~~~\mbox{and}~~~t^\kappa \pp\left(\witi
W_n>t|\calf_0\right)<L_n,~\forall t>0, \feqn and \beqn
\label{plus-plus} \lim_{t \ra \infty} t^\kappa \pp\left(\witi
W_n>t|\calf_0\right)=\witi K_n(\omega).\feqn
\end{proposition}
\begin{remark}
\label{rem28} \item [(i)]  The proof in \cite{kks} of the i.i.d. analogue of 
Proposition \ref{st39} works nearly verbatim with Conditions B and
$\mbox{C}_\kappa$ compensating for the lack of independence of
$\omega.$ Nevertheless, since the proof is rather long and technical,
its detailed modification is included in Appendix
\ref{proof-kks6}. \item [(ii)] The proposition remains valid with the
random variables $\witi W_n$ replaced by the variables $\widehat
W_n=\sum_{j=1}^n W_n.$ The proof is essentially the same, the only
(obvious) difference being that Proposition \ref{moment} can be
applied directly instead of
\eqref{moment3}.  \item [(iii)] Just as with Corollary
\ref{momentc} and Lemma \ref{st17} (a), Proposition \ref{st39}
implies the corresponding uniform estimates for the tails
$\pp\left(\witi W_{m+n} -\witi W_m>t|\calf_{\bar \nu_m} \right)$
as well, for every $m \geq 1.$
\end{remark}
By the bounded convergence theorem, \eqref{plus} and
\eqref{plus-plus} yield \beqn \label{thus} \lim_{t \ra
\infty}t^\kappa \pp\bigl(\witi W_n>t\bigr)=E_P \bigl(\witi
K_n\bigr) \in (0,\infty).\feqn  Note that if either \eqref{finite}
or \eqref{lem-lem} holds with $m=1,$ the random variables $\ol
W_n$ are independent, and the limit laws for their partial sums
follow from the standard i.i.d. limit laws \cite{kolmog,
samor-taqqu}. More generally, we have:
\begin{proposition} \label{w-gvul} Let Assumption \ref{measure1} hold.
\item[(a)] Assume that $\kappa \neq 1.$ Let
$B_n=n^{1/\kappa}$ if $\kappa \in (0,2),$ $B_n=(n \log n)^{1/2}$
if $\kappa=2$,  and $A_n=0$ if $\kappa \in (0,1),$ $A_n= n\ee (\ol
W_2)$ if $\kappa \in (1,2].$ Then, $\left(\witi
W_n-A_n\right)/B_n$ converges in distribution to a stable law of
the form \eqref{kappa-law}. \item[(b)] Assume that $\kappa=1.$
Then, there exist a sequence $\witi D(n) \sim \log n$ and a
positive constant $\tilde c_0$ such that the law of
$\fracd{1}{n}\left(\witi W_n-\tilde c_0 n \witi D(n)\right)$
converges to a stable law of the form \eqref{kappa-law}.
\end{proposition}
\begin{proof}
The random variables $\ol W_n$ are identically distributed and
one-dependent for $n \geq 2$ (see the summary after
\eqref{di-blocks}, and note that we start from $n=2$ because the
slightly different law of $\ol W_1$). Clearly, it is sufficient to
show that the appropriately normalized and centered sums
$S_n=\sum_{j=2}^n \ol W_j$ converge to a stable law of the form
\eqref{kappa-law}. For $\kappa<2,$ apply \cite[Corollary
5.7]{kobus-gp}, noting that the uniform estimates of Proposition
\ref{st39} imply that \beq \forall~\veps>0,\,~ \forall~j\geq
3,~~~~~~ n\pp \left(\ol W_2 \geq \veps n^{1/\kappa}, \ol W_j \geq
\veps n^{1/\kappa} \right)~ \to_{n \to \infty} 0,\feq which is the
tail condition needed to apply Corollary 5.7 of Kobus
\cite{kobus-gp}.

In the case $\kappa=2,$ we note first that $\ol W_2$ and $\ol
W_2+\ol W_3$  both belong by Proposition \ref{st39} to the domain
of attraction of a normal distribution. We seek to apply the limit
theorem in \cite[p. 328]{szewczak}, for which we need to check
that $S_2=\ol W_2$ and $S_3=\ol W_2+\ol W_3$ have different
parameters $b_i=\lim_{n \ra \infty} t^\kappa\pp\left(S_i>t\right),
i=2,3.$ But, \beqn \label{tosefet3} \nonumber b_3&=& \lim_{t \ra
\infty} t^\kappa \pp\left(\ol W_2+\ol W_3>t\right) \geq \lim_{t
\ra \infty} t^\kappa \pp\left( \ol W_2<t,\ol W_3>t\right)+\lim_{t
\ra \infty} t^\kappa \pp\left( \ol W_3<t,\ol W_2>t\right)
\nonumber \\&=& \lim_{t \ra \infty} t^\kappa \pp\left( \ol
W_3>t|\ol W_2<t\right)\pp\left(\ol W_2<t\right)+\lim_{t \ra
\infty} t^\kappa \pp\left(\ol W_2>t\right)\pp\left( \ol W_3<t|\ol
W_2>t\right) \nonumber \\&\geq& J_1+b_2>b_2,\feqn where $J_1$ is
the constant appearing in \eqref{plus}, and we used the uniform
exponential estimates of Proposition \ref{st39} and the fact that
$\pp\left( \ol W_3<t|\ol W_2>t\right) \to_{t \to \infty}1$ which
is also implied by these estimates, as can be seen by conditioning
on the environment to the right of $-\bar\nu_2$. Here and in the
remainder of the proof, any reference to Proposition \ref{st39}
actually includes Remark~\ref{rem28}~{\em (iii)}. We have \beqn
\label{tosefet1}  \lim_{t \to \infty}\pp\left( \ol W_3<t|\ol
W_2>t\right)= \lim_{t \ra \infty} \ee\left(\pp\left( \ol W_3<t|
\calf_{\bar \nu_2} \right)|\ol W_2>t\right).\feqn By Proposition
\ref{st39}, \beq \pp\left( \ol W_3<t| \calf_{\bar \nu_2} \right)
\geq 1-L_1 t^{-\kappa},~~~P-\as, \feq implying that the limit in
\eqref{tosefet1} exists and is equal to 1. Therefore, by
\eqref{tosefet3} and since we know a-priori from \eqref{plus-plus}
that $b_3=\lim_{t \to \infty} t^\kappa \pp(\ol W_2+\ol W_3>t)$ is
well-defined,  the following limit exists and can be bounded below
by using \eqref{plus}: \beq \lim_{n \ra \infty} t^\kappa\pp\left(
\ol W_3>t|\ol W_2<t\right)= \lim_{t \ra \infty}t^\kappa
\ee\left(\pp\left( \ol W_3>t| \calf_{\bar \nu_2} \right)|\ol W_2<t
\right) \geq J_1.\feq This completes the proof of the proposition.
\end{proof}
\paragraph{Completion of the proof of Proposition  \ref{markov-tau}.}
The limit laws for $T_n$ announced in Proposition \ref{markov-tau} are
obtained from stable laws for partial sums of $\ol W_n$ in the same
way as in \cite{kks}, by a standard argument using Lemma
\ref{st17}. To illustrate the argument we consider here the case
$\kappa=2,$ omitting the proof for $\kappa \in (0,2).$ Let
$\zeta(n)=\max\{i: \bar  \nu_i<n\}$ and
$\varsigma(n)=[n/\mu-C\sqrt{n}]$ for a constant $C>0.$ Using part
{\em (c)} of Lemma \ref{st17}, we obtain, with $\mu=\ee(\bar
\nu_2-\bar\nu_1)$, \beq \liminf_{n\ra\infty} \pp\left(\zeta(n)\geq
n/\mu-C\sqrt{n}\right) &\geq& \lim_{n \ra \infty} \pp\left(\bar
\nu_{\varsigma(n)} \leq n\right)
\\&=& \lim_{n \ra \infty} \pp\left(\frac{\bar \nu_{\varsigma(n)}-
\varsigma(n) \mu}{\sigma \sqrt{\varsigma(n)}} \leq
\frac{n-\varsigma(n)\mu}{\sigma
\sqrt{\varsigma(n)}}\right)=\call_{2,\frac{\sigma}{\sqrt{2}}}
\left(C \mu^{3/2}\right).\feq Hence, for all $\epsilon>0$ and some
$C=C(\veps)>0$ and all $n>N_2(\veps),$ $\pp\left(\zeta(n) \leq
n/\mu-C\sqrt{n}\right) \leq \veps .$ It follows, letting $a= \ee
(\ol W_2),$ that for any $n$ large enough, \beq &&\pp
\left(\frac{\sum_{i=1}^{n} Z_i -na/\mu}{\sqrt{n \log n}} \leq
x\right) \leq \pp \left(\sum\limits_{i=1}^{\zeta(n)} \ol W_i \leq
x \sqrt{n \log n}+na/\mu\right) \\&&\leq\pp
\left(\sum\limits_{i=1}^{[n/\mu -C\sqrt{n})]} \ol W_i \leq x
\sqrt{n \log n}+na/\mu \right)+\veps \ra \call_{2,\tilde b}
\left(x \sqrt{\mu}\right)+\veps, \feq where $\call_{2,\tilde b}$
is the limiting law for sums of $\ol W_n.$ Similarly, \beq &&\pp
\left(\frac{\sum_{i=1}^{n} Z_i-na/\mu}{\sqrt{n \log n}} \leq
x\right)  \geq \pp \left(\sum\limits_{k=1}^{\zeta(n)+1} \ol W_k
\leq x \sqrt{n \log n} +na/\mu\right) \\&&\geq\pp
\left(\sum\limits_{k<n/\mu +C\sqrt{n}} \ol W_k \leq x \sqrt{n \log
n}+an/\mu\right)-\veps\ra \call_{2,\tilde b} \left(x \sqrt{\mu}
\right)-\veps. \feq Since $\veps$ was arbitrary, Proposition
\ref{markov-tau} now follows from the limit laws  for partial sums
of $Z_n$ by \eqref{ti-zed}--\eqref{zed-sums}. Since the law
defined by \eqref{kappa-law} has expectation zero,
$\mbox{v}_P=a/\mu=\ee(\tau_1),$ where $\tau_1$ is defined by
\eqref{ti-tau}. \qed \\ As shown in the Introduction this
completes the proof of Theorem \ref{main-markov}.
\subsection{Tails of distribution of the random variable $R$}
\label{proof-ar} The aim of this subsection is to prove that
Condition $\mbox{C}_\kappa$ holds for some $\kappa>0.$ Proposition
\ref{sb} below extends the following  theorem, valid in the i.i.d.
setup, to some Markov-dependent variables.
\begin{theorem}[Kesten]\cite[Theorem 5]{kesten-randeq} \label{goldenkest}
Let $(Q_n,M_n),~n \in \nn,$ be independent copies of a $\rr^2$-valued
random vector $(Q,M),$ satisfying the following conditions:
\begin{itemize} \item[(i)] $P(M>0)=1$ and $P(Q>0)=1.$
\item[(ii)] For some $\kappa>0,$ $E\left(M^\kappa\right)=1,$
$E\left(M^\kappa \log^+ M\right)<\infty,$ and $E(Q^\kappa)<\infty.$
\item[(iii)] The law of $\log M$ is non-lattice (its support is not 
contained in any proper sublattice of $\rr$) and
$P(Q=(1-M)c)<1,~\forall c \in \rr.$
\end{itemize}
Then there exists a constant $\widehat K>0$ such that \beqn
\label{gvul-kova} \lim_{t \ra \infty} t^\kappa P(\widehat R \geq
t)=\widehat K, \feqn where $\widehat
R:=Q_1+M_1(Q_2+M_2(Q_3\ldots)).$
\end{theorem}
We have:
\begin{proposition}\label{sb}
Let Assumption \ref{measure1} hold.  Then Condition
$\mbox{C}_{\kappa}$ is satisfied for the $\kappa>0$ defined by
\eqref{kappa239}.
\end{proposition}
\begin{proof}
If either \eqref{finite} or \eqref{lem-lem} with $m=1$ hold, this
proposition can be deduced rather directly from Kesten's theorem.
It will be convenient to give a separate proof for the case where
the state space $\cals$ is finite, i.e. under assumption
\eqref{finite}. \item[{\bf Assume first that \eqref{finite}
holds}.] Then, it is sufficient to show that \beq K_x:=\lim_{t \ra
\infty} t^\kappa P_x(R
>t) \in (0,\infty)\feq exists for all $x \in \cals.$ For $n \geq 0,$ let \beqn
\label{qm} Q_n=1+\odin_{\{N_{n+1} \geq N_n+2\}}
\sum_{i=N_n}^{N_{n+1}-2} \prod_{j=N_n}^i \rho_{-j}~~~\mbox{and}~~~
M_n=\prod_{i=N_n}^{N_{n+1}-1} \rho_{-i}. \feqn Then, $(M_n,Q_n)_{n
\geq 1}$ is an i.i.d. sequence, and
$R=Q_0+M_0(Q_1+M_1(Q_2+\ldots))$. First, we will show that
Kesten's theorem is applicable to this sequence, that is the
following limit \beqn \label{rkova} \widehat K:=\lim_{t \ra
\infty} t^\kappa P_x(\widehat R >t) \in (0,\infty) \feqn exists,
where \beqn \label{rkova1} \widehat
R=Q_1+M_1(Q_2+M_2(Q_3\ldots))\,, R=Q_0+M_0\widehat R.\feqn Let
$f_\kappa$ be a strictly positive Perron-Frobenius eigenvector of
the matrix $H_\kappa(x,y):=~$ $H(x,y)\rho(y)^\kappa.$ By virtue of
\eqref{hbeta} and Condition B, it corresponds to the eigenvalue
$1.$ Recall now the definitions of the state $x^*$ and the matrix
$\Theta$ from \eqref{kernels1}. By Lemma \ref{urg1}, the
Perron-Frobenius eigenvalue (the spectral radius) of the matrix
$\Theta_\kappa(x,y)=\Theta(x,y)\rho(y)^\kappa$ is strictly less
than one. So, the vector $f_\kappa$ normalized by the condition
$f_\kappa(x^*) \rho^\kappa (x^*)=1$ is the unique positive vector
in $\rr^{|\cals|}$ solving the equation $(I-\Theta_\kappa)f=s,$
where $s(x):=H(x,x^*).$ Hence (this is a very particular case of
the results of \cite{athreya-ney1} and \cite[Theorem
5.1]{nummelin} ) \beqn \label{st} f_\kappa(x)= \rho(x)^{-\kappa}
E_x\left(\prod_{i=0}^{N_1-1}\rho_{-i}^\kappa\right)=
\sum_{n=0}^\infty \Theta_\kappa^n s(x)\,,\feqn and \beqn
\label{st11} E_{x^*} \left
(\prodl_{i=0}^{N_1-1}\rho_{-i}^\kappa\right)=E_P(M_1^\kappa)=1.
\feqn The second equality in \eqref{st} follows since the chain
$(x_i)$ evolves according to the kernel $\Theta$ until $N_1$ (see
\eqref{between}), while \eqref{st11} follows from the
normalization condition $f_\kappa(x^*) \rho^\kappa(x^*)=1.$

It is not hard to check that assuming \eqref{finite}, condition {\em
(A4)} is equivalent to the fact that $\log M_1$ is non-lattice, and
that $P(Q_1=(1-M_1)c)<1$ for any $c \in \rr$\ \ (since clearly
$P(M_1>1)>0$), as required to apply Theorem~\ref{goldenkest}. In order
to prove \eqref{rkova}, it remains to show that
$E_P(Q_1^\kappa)<\infty$ and $E_P\left(M_1^\kappa \log^+
M_1\right)<\infty.$ Thus, it is sufficient to prove that there exists
$\beta>\kappa$ such that \beqn \label{233} E_x\bigl(Q_0^\beta
\bigr)~\mbox{is a bounded function of $x$}.
\feqn Since for any $n \in \nn$ and positive numbers
$\{a_i\}_{i=1}^n$ we have \beq (a_1+a_2+\ldots a_n)^\beta \leq
n^\beta (a_1^\beta+a_2^\beta+\ldots a_n^\beta), \feq we obtain for
any $\beta>0$ and $x \in \cals:$ \beqn \label{s33}
E_x\bigl((Q_0-1)^\beta \bigr)& =&E_x \left(\sum_{n=2}^\infty
\sum_{i=1}^{n-1} \prod_{j=0}^{i-1} \rho_{-j}
\one{N_1=n}\right)^\beta =\sum_{n=2}^\infty E_x \left(
\sum_{i=1}^{n-1} \prod_{j=0}^{i-1} \rho_{-j}
\one{N_1=n}\right)^\beta \nonumber \\& \leq &\sum_{n=2}^\infty
(n-1)^\beta \sum_{i=1}^{n-1} E_x \left(\prod_{j=0}^{i-1} \rho_{-j}
^\beta \one{N_1 \geq n}\right). \feqn But $E_x
\left(\prod_{j=0}^{i-1} \rho_{-j} ^\beta \one{N_1 \geq
n}\right)=\rho(x)^\beta \Theta^{n-i}\Theta_\beta^{i-1} \odin,$
where $\Theta_\beta(x,y):=\Theta(x,y)\rho(y)^\beta.$ Since the
spectral radius of the matrices $\Theta_\kappa$ and $\Theta$ are
strictly less than one, it follows from \eqref{s33} that
\eqref{233} holds for some $\beta>\kappa.$ This yields
\eqref{rkova}.

By \eqref{rkova} and the bounded convergence theorem, and since
the random variables $M_0$ and $\widehat R$ are independent under
the measure $P_x,$ the following limit exists: \beq K_x:= \lim_{t
\ra \infty} t^\kappa P_x(M_0 \widehat R>t)=\widehat K
E_x(M_0^\kappa) \in (0,\infty).\feq Fix any $\alpha \in
\bigl(\frac{\kappa}{\beta}, 1\bigr).$ It follows from
\eqref{rkova} and \eqref{233} that for all $t>1,$ \beq t^\kappa
P_x( R>t)&\leq& t^\kappa P_x\bigl(Q_0+M_0\widehat
R>t,~Q_0<t^\alpha\bigr)+t^\kappa P_x (Q_0 \geq t^\alpha)
\nonumber\\ &\leq& t^\kappa P_x\bigl(M_0\widehat
R>t-t^\alpha\bigr)+\fracd{t^\kappa}{t^{\alpha \beta}} E_x
(Q_0^\beta), \feq and \beq && t^\kappa P_x( R>t)=t^\kappa
P_x\bigl(Q_0+M_0\widehat R>t\bigr) \geq t^\kappa
P_x\bigl(M_0\widehat R>t\bigr).\feq We conclude, by taking the
limit in the above inequalities as $t \ra \infty,$ that \beq
\lim_{t \ra \infty} t^\kappa P_x(R>t)=\lim_{t \ra \infty} t^\kappa
P_x(M_0 \widehat R>t)=K_x,\feq completing the proof of the
proposition in the case \eqref{finite}. \item[{\bf Assume now that
\eqref{lem-lem} holds}.] First, we will prove that \eqref{gvul}
holds for some function $K(\omega)$ and constant $\widehat K.$ We
follow Goldie's proof \cite{goldie} of Kesten's Theorem
\ref{goldenkest}. Let \beq \eta(x):=\log \rho(x),\feq \beq
\Pi_0=1,~\Pi_n=\prod_{k=0}^{n-1} \rho_{-k},~~~~n \geq 1, \feq \beq
\eta_n=\log \rho_{-n},~V_n=\log \Pi_n~~~n \geq 0, \feq \beqn
\label{kohav} R=R^0=\sum_{n=0}^\infty
\Pi_n,~~~R_0=0,~R_n=\sum_{k=0}^{n-1} \Pi_k,~~~
R^n=(R-R_n)/\Pi_n,~~~ n \geq 1.  \feqn Following Goldie
\cite{goldie}, we write for any numbers $n \in \nn,~t \in \rr,$
and any point $z \in \cals,$ \beq P_z(R
>e^t)&=&\sum_{k=1}^n [ P_z(e^{V_{k-1}}
R^{k-1}>e^t)-P_z(e^{V_k} R^k>e^t)]+ P_z(e^{V_n} R^n>e^t).\feq We
have, by using the identity $R^{k-1}=1+\rho_{-k+1}R^k,$ \beq
&&P_z(e^{V_{k-1}} R^{k-1}>e^t)-P_z(e^{V_k} R^k>e^t)=
\\&& \int_\rr \int_\cals [P(R^{k-1}>e^{t-u}|x_{k-1}=x)
-P(\rho_{_{-k+1}}R^k>e^{t-u}|x_{k-1}=x)] P_z(V_{k-1} \in
du,x_{k-1} \in dx)\\&& = \int_\rr \int_\cals
[P_x(R>e^{t-u})-P_x(R-1>e^{t-u})] P_z(V_{k-1} \in du,x_{k-1} \in
dx). \feq Thus, letting $\delta_n(z,t)= e^{\kappa t} P_z(e^{V_n}
R^n>e^t)$ and $f(x,t)=e^{\kappa t}[P_x (R>e^t)-P_x (R-1>e^t)],$
\beqn \label{for-sum} r_z(t):=e^{\kappa t}P_z(R
>e^t)=\sum_{k=0}^{n-1} \int_\rr \int_\cals f(x,t-u) e^{\kappa u}
P_z(V_k \in du,x_k \in dx)+\delta_n(z,t). \feqn By Lemma \ref{urg} and
\eqref{hbeta}, there exists a positive measurable function $h(x):
\cals \ra \rr$ bounded away from zero and infinity such that: \beq
h(x)=\int_\cals H(x,dy)\rho^\kappa (y) h(y). \feq This implies, by
\cite[Theorem 5.2]{nummelin}, that there is a probability measure
$\pi_\kappa$ invariant for the kernel
$H_\kappa(x,dy)=H(x,dy)\rho^\kappa(y)$, namely (since
$r_{_{H_\kappa}}=1$ by \eqref{kappa239} and \eqref{hbeta}) \beqn
\label{eq} \int_\cals H_\kappa(x,A)
\pi_\kappa(dx)=\pi_\kappa(A),~~~ \forall~A \in \calt. \feqn The
measure $\pi_h(dx)=h(x)\pi_\kappa(dx)$ is a finite invariant
measure for the kernel \beq \witi H(x,dy):=\fracd{1}{h(x)}
H_\kappa(x,dy)h(y). \feq The measure $\pi_\kappa$ and hence
$\pi_h$ are equivalent to the original stationary distribution
$\pi.$ Indeed, by \eqref{eq}, \beq \int_\cals H_\kappa^m(x,A)
\pi_\kappa(dx)=\pi_\kappa(A),~~~ \forall~A \in \calt. \feq Hence,
by \eqref{lem-lem} and the ellipticity condition {\em (A2)},
$c_r^{-1}c_\rho^{-m}\pi_\kappa(A) \leq \pi(A) \leq c_r c_\rho^m
\pi_\kappa(A),$ where the constant $c_\rho$ is defined in \eqref{crho}.

Let $\witi P$ be the probability measure under which the Markov chain
$(x_k)_{k \geq 0}$ is stationary and governed by the transition
probability measure $\witi H(x,A).$ As usual we denote the conditional
probabilities $\witi P(\cdot|x_0=x)$ by $\witi P_x(\cdot).$ Then,
\beq r_z(t)=\sum_{k=0}^{n-1}
\int_\rr \int_\cals f(x,t-u) \frac{\rho^\kappa(z) h(z)}{\rho^\kappa
(x) h(x)} \witi P_z(V_k \in du,x_k \in dx) +\delta_n(z,t). \feq Since
$P-\as,$ $\Pi_n R^n \ra 0$ as $n$ goes to infinity, $P\bigl(\lim_{n\ra
\infty} \delta_n(z,t)=0\bigr)=1,$ for any fixed $t>0$ and $z \in
\cals.$  Therefore, $P-\as,$ \beq r_z(t):=e^{\kappa t}P_z(R
>e^t)=\sum_{k=0}^\infty \int_\rr \int_\cals f(x,t-u)
\frac{\rho^\kappa (z) h(z)} {\rho^\kappa (x) h(x)} \witi P_z(V_k
\in du,x_k \in dx).\feq We will use the following Tauberian lemma
: \begin{lemma}\cite[Lemma 9.3]{goldie} \label{igol}Let $R$ be a
random variable defined on a probability space $(\Omega,\calf,
P).$ Assume that for some constants $\kappa,K \in (0,\infty),$
$\int_0^t u^\kappa P(R>u)du \sim Kt$ as $t \ra \infty.$ Then
$t^\kappa P(R>t) \sim K.$
\end{lemma} It follows from Lemma \ref{igol} that in
order to prove \eqref{gvul}, it is sufficient to show that $P-\as$
there exists \beq \lim_{t \ra \infty} \check{r}_z(t) \in
(0,\infty), \feq where the smoothing transform $\check{q}$ is
defined, for a measurable function $q: \rr \ra \rr$ bounded on
$(-\infty,t]$ for all $t,$ by \beq \check{q}(t):=\int_{-\infty}^t
e^{-(t-u)}q(u)du. \feq Let \beqn \label{gi} g(x,t)&:=&
\frac{1}{e^{\kappa \eta(x)} h(x)} \int_{-\infty}^t
e^{-(t-u)}f(x,u)du \nonumber \\ &=& \frac{1}{e^{\kappa \eta(x)}
h(x)} \int_{-\infty}^t e^{-(t-u)} e^{\kappa u}[P_x (R>e^u)-P_x
(R-1>e^u)]du \nonumber \\&=& \frac{e^{-t}}{e^{\kappa \eta(x)}
h(x)}\int_0^{e^t} v^\kappa[P_x (R>v)-P_x (R-1>v)]dv.  \feqn Then,
since $\check{r}_z(t)= h(z) \rho^\kappa (z) \sum_{k=0}^\infty
\witi E_z(g(x_k,t-V_k)),$ it is sufficient to show that for any $z
\in \cals,$ \beqn \label{hp} \lim_{t \ra \infty} \sum_{k=0}^\infty
\witi E_z(g(x_k,t-V_k)) \feqn exists and belongs to $(0,\infty).$
So, our goal now is to prove \eqref{hp}.

Toward this end, note first that the kernel $\witi H$ satisfies
condition \eqref{lem-lem} and hence the chain $(x_n)$ is ergodic under
the measure $\witi P.$ Further, the random walk $V_n=\sum_{j=0}^{n-1}
\eta_j$ has a positive drift under the measure $\witi P_x.$ Indeed,
similarly to \cite{goldie} and
\cite{kesten-randeq}, we obtain for some $c>0$ and any $\gamma>0,$
\beq \witi P_x\left(e^{V_n} \leq e^{-\gamma n^{1/4}}\right)&=&
\fracd{e^{-\kappa \eta(x)}}{h(x)} E_x\left(e^{\kappa V_n}
h(x_{n-1});e^{V_n} \leq e^{-\gamma n^{1/4}}\right) \leq c
E_x\left(e^{\kappa V_n};e^{V_n} \leq e^{-\gamma n^{1/4}}\right)
\\&\leq& c e^{-\kappa \gamma n^{1/4}}.\feq Thus, $\lim_{n \ra \infty}
\witi P_x\left(V_n \leq -\gamma n^{1/4}\right)=0,$ implying $\witi
E_{\pi_h} (\eta_0)>0$ by the central limit theorem for bounded
additive functionals of Doeblin recurrent Markov chains (see
e.g. \cite[p. 134]{nummelin}).

The limit in \eqref{hp} follows from the version of the Markov renewal
theorem as given in \cite[Theorem 1]{alsmeyer} (see also
\cite{mac-dan, kesten-rem}) when applied to the Markov chain
$(x_{n+1},\rho_{-n})$, provided that we are able to show that the
following holds: \beqn \label{al1} g(x,\cdot)~\mbox{is a continuous
function for}~\pi_h-\mbox{almost all}~x \in \cals, \feqn and \beqn
\label{al5} \int_\cals \sum_{n
\in \zz} \sup_{n \delta \leq t< (n+1) \delta } \{|g(x,t)|\}
\pi_h(dx)<\infty~~~\mbox{for some}~\delta>0. \feqn

The assertion \eqref{al1} follows from the continuity of $\int_0^{e^t}
v^\kappa[P_x (R>v)-P_x (R-1>v)]dv$ in $t$ for every $x \in \cals.$ For
some $M>0$ and any $\veps \in (0,1),$ we get from \eqref{gi}: \beq
&&g(x,t) \leq M e^{-t}\int_0^{e^t} v^\kappa[P_x (R>v)-P_x (R-1>v)]dv
\\ &\leq& M e^{-\veps t}\int_0^{e^t} v^{\kappa-1+\veps} [P_x (R>v)-P_x
(R-1>v)]dv \leq \frac{M}{\kappa} e^{-\veps t} E_x [(R)^{\kappa+\veps}-
(R-1)^{\kappa+\veps}], \feq where the last inequality follows from
\cite[Lemma 9.4]{goldie}.  Since for any $\gamma>0$ and $R>1,$
$(R)^\gamma- (R-1)^\gamma \leq
\max\{1,\gamma\} (R)^{\max\{1,\gamma\}-1},$ we obtain by Condition
B that \beq E_x [(R)^{\kappa+\veps}- (R-1)^{\kappa+\veps}] \leq L,
\feq for some constant $L>0$ independent of $x,$ yielding
\eqref{al5} and consequently \eqref{gvul}. In fact we
have shown that the following limit exists $\pi-\as :$ \beqn
\label{27x} \lim_{t \ra \infty} t^\kappa P_x (R>t)=K_1(x) \in
(0,\infty). \feqn We now turn to the proof of \eqref{bound}. Fix
any point $x^* \in \cals$ for which \eqref{27x} holds. Using
\eqref{crho} and \eqref{lem-lem}, we obtain for any $x \in \cals$
and $t >0 :$ \beq P_x (R>t) \geq P_x (c_\rho^{-m} R^m >t)
=\int_\cals H^m (x,dz) P_z (c_\rho^{-m} R
>t) \geq c_r^{-2} P_{x^*} (c_\rho^{-m} R
>t), \feq and \beq P_x (R>t) &\leq& P_x (m c_\rho^m+c_\rho^m R^m
>t) =\int_\cals H^m (x,dz) P_z (m c_\rho^m+ c_\rho^m R >t) \\
&\leq& c_r^2 P_{x^*} (m c_\rho^m+c_\rho^m R >t). \feq Thus,
\eqref{bound} follows from \eqref{27x}.
\end{proof}
\begin{remark}
It should be mentioned that essentially the same proof leads to
similar tail estimates for random variables of the form
$R=\sum_{n=0}^\infty Q_n \prod_{j=0}^{n-1}M_j$ with a more general
type of Markov-dependent coefficients $(Q_n,M_n)$ (e.g. $Q_n$ need
not be deterministic and $M_n$ need not be a.s. positive). This
general result (under somewhat milder assumptions than those  assumed
in this paper, namely allowing for periodic Markov chains  while
relaxing the uniform bound \eqref{lem-lem} on the kernels
$H(x,\cdot)$) can be found in \cite{rec}. While preparing the  final
version of the article for publication, we were kindly  informed by
J.~Bremont of the reference \cite{saporta} where, by  using different
methods, a result similar to Proposition \ref{sb}  is obtained for
Markov chains in a finite state space.
\end{remark}
\section{Summary and Final Remarks}
\label{conlus}
 We have dealt with the random walk $(X_n)_{n\ge 0}$ in a random
 environment $\omega\!\in\![0,1]^{\mathbb{Z}}$, associating with it an
 auxiliary Galton--Watson process $(Z_k)_{k\ge 0}$ with one
 immigrant at each instant and random branching mechanism\
 ${\rm Geom}(\omega_{-k})$.\\
\centerline{}
\\
\noindent Without stating it explicitly the following theorem has in 
fact been proved.  Let $(\calg_n)_{n \in \nn}$ be an augmentation of
$(\sigma(\omega_j:j>-n))_{n \in \nn}$ which generates the original
quenched law, namely $\pp(\cdot|\calg_n)=P_\omega(\cdot),$ $\as.$
Accordingly, let $(\ol \calf_n=\sigma(Z_0,Z_1,\ldots,Z_n) \vee
\calg_n)_{n \in \nn}$ be $(\calf_n)$'s induced augmentation.
\begin{theorem}
\label{conc-thm} Assume the environment $\omega$ satisfies Conditions
B and C$_\kappa$ (for the $\kappa>0$ involved in Condition B)
introduced in Section~\ref{proofs}. Furthermore, assume the existence
of an increasing sequence of stopping times $\eta_m$, with $\eta_0=0$,
with respect to the filtration $(\ol \calf_n)_{n\ge 0}$ (defined
in~\eqref{f-zero} and~\eqref{sigmag}) for which \item[i)] the LLN and
CLT hold: there exist\ \ $\mu>0$\ \ and\ \ $\sigma\!\in\!\mathbb{R}$\
\ such that \[ \frac{\eta_m}{m}\longrightarrow\mu\ \ \ {\rm a.s.}\ \ \
\ \ {\rm and} \ \ \ \ \ \frac{\eta_m-m\mu}{\sqrt{m}}\ \ \stackrel{\cal
D}{\longrightarrow}\ \ N(0,\sigma^2) \] \item[ii)]\ for some $b>0,
\qquad \fracd{1}{B_m}\left(\sum_{_{i=1}}^{^{\eta_{_m}}}Z_i\ -\
A_m\right)\ \ \stackrel{\cal D}{\longrightarrow}\ \ \call_{\kappa,b}
\quad ({\rm defined\ in}~\eqref{kappa-law})$\\
\centerline{}
\\
    \noindent
    where $A_m\left\{\begin{array}{ll}
                 =0               & \kappa\in(0,1) \\
                 \sim c_1 m\log m & \kappa=1       \\
                 =c_2 m           & \kappa\in(1,2] \\
                \end{array} \right.$
  {\rm and}
          $B_m=\left\{\begin{array}{ll}
                 m^{\frac{1}{\kappa}}     & \kappa\in(0,2) \\
                 (m\log m)^{\frac{1}{2}} & \kappa=2      \\
                \end{array}  \right.$\\
\centerline{}
\\
   \noindent
   for suitable positive constants $c_1,c_2$.
   \\
 \centerline{}
\\
 \noindent
Then the random walk $X_n$ satisfies a stable limit law in the
sense that the conclusions $(i)$--$(iv)$ of
Theorem~\ref{main-markov} hold.
 \end{theorem}
 \bigskip

In the Markov setup of this paper, and under Assumption
\ref{measure1}, we have shown (see Lemma~\ref{st17} and
Proposition~\ref{w-gvul}) that the environment $\omega$ indeed
satisfies the conditions of Theorem~\ref{conc-thm} (with respect to
the stopping times $\eta_n=\bar \nu_n$), thus obtaining the stable
limit laws in this case.

It is easy to see that Theorem \ref{main-markov} can be extended
for instance to the following setup of hidden Markov models. Let
$(x_n)_{n \in \zz}$ be a Markov chain defined on a measurable
space $(\cals,\calt)$ that satisfies {\em (A1)} and {\em (A2)} in
Assumption~\ref{measure1}. Assume that in the product space
$(\cals \times \Omega,\calb \times \calf),$ \beqn \label{hmc}
P\bigl((x_n,\omega_{-n}) \in A \times B| x_{n-1}=x,
\sigma((x_i,\omega_{-i}): i \leq n-1)\bigr)={\mathbb H}(x,A \times
B)\feqn for all $n \in \zz, A \in \calt, B \in \calf, x \in
\cals,$ where ${\mathbb H}$ is a stochastic kernel on $(\cals,
\calt \times \calb).$ Note that the Markov chain
$(x_n,\omega_{-n})$ might not satisfy Assumption \ref{measure1},
so that Theorem \ref{main-markov} cannot be applied directly.

Let $Q(x,y,B)=P(\omega_{-n} \in B|x_{n-1}=x,x_n=y).$ Then, similarly
to \eqref{hbeta}, 
\beq E_x \left(\prod_{i=0}^{n-1} \rho_{-i}^\beta \right)=
\rho_0^\beta H^{n-1}_\beta \odin(x),~~~\beta \geq 0, 
\feq where the kernel $H_\beta(x,\cdot)$ is
now defined on $(\cals,\calt)$ by \beqn \label{kbeta}
H_\beta(x,dy)=H(x,dy) \int_\Omega Q(x,y,dz) \rho^\beta(z). \feqn From
the ellipticity condition {\em (A2)} it follows that $\rho_0
\in (c_\rho^{-1},c_\rho)$ for some constant $c_\rho>0,$ and we
obtain that $ C_\beta^{-1} H(x,dy) \leq H_\beta(x,dy) \leq C_\beta
H(x,dy),$ for a suitable constant $C_\beta>0.$ Thus, Lemma
\ref{urg} is in force for the kernel $H_\beta$ defined by
\eqref{kbeta}.

The non-arithmetic condition needed to apply a Markov renewal
theorem (\cite[Theorem 1]{alsmeyer}) now takes the form:
\begin{definition}
\label{shurd} \cite{shur,alsmeyer} The process $\log \rho_{-n}$ is
called $\alpha$-arithmetic if $\alpha \geq 0$ is the maximal
number for which there exists a measurable function $\gamma :\cals
\ra [0,\alpha)$ such that \beq P\left( \log \rho_0 \in
\gamma(x)-\gamma(y)+\alpha \cdot
\zz|x_{-1}=x,x_0=y\right)=1,~~~P-\as \feq The process is called
non-arithmetic if no such $\alpha$ exists.
\end{definition}
We have:
\begin{theorem}
\label{th-first} Assume that the underlying model
$\omega_n=\omega(x_n)$ in Theorem \ref{main-markov} is replaced by
\eqref{hmc} an that Assumption \ref{measure1} holds. Then the 
conclusions $(i)$--$(iv)$ of Theorem~\ref{main-markov} remain valid.
\end{theorem}
The proof is the same as that of Theorem \ref{main-markov} by
using the regeneration times $\bar \nu_n$ defined in
\eqref{nubar}. The only exception is that in the definition of
$f(x,t)$ (a line before \eqref{for-sum}) we would condition on
$x_{-1}$ rather than on $x_0.$ Correspondingly, in the definition
of $r_\lambda$ (cf. \eqref{for-sum}), the integration would be
with respect to the measure $P_\lambda(V_k \in dv, x_{k-1} \in
dx).$
\begin{appendix}
\section*{Appendix}
Recall $\f_0=\sigma(\omega_k: k>0).$ For brevity, we denote the
conditional probabilities $P(\cdot|\f_0)$ and $\pp(\cdot|\f_0)$ by
$P_+$ and $\pp_+=P_+ \otimes P_\omega$ respectively. We usually do not
indicate the argument $\omega$ of these functions meaning that the
inequalities below hold $P-\as$ We denote by $\theta$ the shift on
$\Omega$, given by $(\theta\omega)_i=\omega_{i+1}.$ For an event $A,$
the notation $I(A)$ is used to denote the indicator function of $A.$
\section{Proof of Proposition \ref{moment}} \label{moment-proof}
The key to the proof is
\begin{lemma}\cite[(2.12)]{kks}
\label{ladder} Suppose that the environment $\omega$ is stationary
and ergodic, and $a_P:=E_P(\log \rho_0)<0.$ Choose any $\gamma \in
(a_P,0)$ and define \beq U_n&=&\sum\limits_{i=0}^{n-1}\left\{\log
\rho_{-i}
-\gamma\right\}~~~(U_0=0), \nonumber \\
\zeta_0&=&0,~~~\zeta_{k+1}=\inf\{n>\zeta_k:U_n \leq U_{\zeta_k}\}.
\feq Then there exist constants $K_5,K_6>0$ such that $P-\as,$
\beq P_\omega(\nu_1>\zeta_k) \leq K_5e^{-K_6k},~~~~k> 0. \feq
\end{lemma}
\begin{remark}
This lemma is proved in \cite{kks} for the special case
$\gamma=a_P/2,$ but an inspection of the proof reveals that
$a_P/2$ can be replaced by any constant between $a_P$ and zero in
the definition of the random walk $U_n.$
\end{remark}
\noindent By virtue of Lemma \ref{ladder}, it is sufficient to
find $\gamma<\in (a_P,0)$ such that for some constants $b>0$ and
$K_7,K_8>0$ \beq P_+(\zeta_k > b k)<K_7e^{-K_8k},~~~~k \geq 0.
\feq Let $\eta(n)=\max\{j: \zeta_j \leq n\}$ and recall
$c_\rho=(1-\epsilon)/\epsilon.$ Since for any $n>0,$ \beq U_n \geq
\sum_{j=1}^{\eta(n)} (U_{\zeta_j}-U_{\zeta_{j-1}}) \geq -\eta(n)(
\gamma+\log c_\rho), \feq for any $k>0,$ the event $\{\zeta_{k+1}
>n \}=\{\eta(n) \leq k\}$ is included in $\{U_n \geq -k \log c_\rho-k
\gamma\}.$ Therefore, for any $\gamma \in (a_P,0)$ and $b \in \nn$
we have \beq P_+(\zeta_{k+1} >bk)  \leq
P_+\left(\sum\limits_{i=0}^{bk-1} \log \rho_{-i} \geq -k\log
c_\rho+k(b-1)\gamma \right).  \feq Let $\gamma=\fracd{1}{2}\cdot
\lim\limits_{\beta \ra +0}\fracd{\Lambda(\beta)}{\beta},$ where
$\Lambda(\beta)$ is as in \eqref{kappa239}, noting that since
$\Lambda(\beta)$ is convex, $\gamma$ is negative by Condition B
and is greater than $a_P$ by Jensen's inequality. Hence, by
Chebyshev's inequality and Condition B, we obtain for any fixed
$b>0$ and $\beta>0$ small enough, \beq \limsup_k \fracd{1}{k} \log
\left[ P_+(\zeta_{k+1}
> k b )\right] &\leq&  \beta \log c_\rho
-(b-1)\gamma\beta +\fracd{3b\gamma \beta}{2}=\beta \left(\log
c_\rho+\gamma+\fracd{b\gamma}{2} \right). \feq Taking $b>-4\log
c_\rho /\gamma$ in the last inequality gives
\begin{equation*}
\limsup_k \fracd{1}{k} \log P_+\bigl(\zeta_{k+1} > k b\bigr) <
\beta (-\log c_\rho+\gamma)<0.
\end{equation*}
This completes the proof of Proposition \ref{moment}.
\section{Proof of Proposition \ref{st39}}
\label{proof-kks6} As mentioned in Remark \ref{rem28} {\em (i)},
this proof will follow the one of \cite[Lemma 6]{kks} very closely, at
times word by word, with the necessary changes made in annealed
arguments to take the dependence of the environment into
account. Quenched arguments, where no changes are needed, will be
skipped.

Throughout we fix a number $\tilde n \in \nn$ and denote $\witi
W:=\witi W_{\tilde n}=\sum_{j=1}^{\tilde n} \ol W_j,$ $\tilde
\nu:=\bar \nu_{\tilde n}.$ Recall the filtration $(\calf_n)_{n \geq 0}$
introduced in \eqref{f-zero} and \eqref{sigmag}, and for all $A>0$
define its stopping time $\varsigma_A=\inf\{n: Z_n>A\}.$ The random
variable $\witi W$ can be represented on the event
$\{\varsigma_A<\tilde \nu\}$ in the following form: \beqn
\label{piruk-matrif} \witi W=\suml_{n=0}^{\varsigma_A-1}Z_n+
S_{\varsigma_A}+\suml_{\varsigma_A \leq n < \tilde \nu} Y_n, \feqn
where \beq Z_{n,k}&=&\mbox{number of progeny alive at time $k$ of the
immigrant who entered at time}~n<k,\\
Y_n&=&\suml_{k>n}Z_{n,k}=\#\{\mbox{progeny of the immigrant at time
$n,$ not including the immigrant}\} \\S_n&=&Z_n+\mbox{total progeny of
the $Z_n$ particles present at}~n. \feq It will turn out that for
large $A,$ the main contribution to $\witi W$ in
\eqref{piruk-matrif} comes from the second term and $\pp_+\bigl(
\witi W \geq t \bigr)\approx \pp_+\bigl(S_{\varsigma_A} \geq
t,\varsigma_A<\tilde \nu \bigr).$ If an environment $\omega$ is
fixed, then $S_{\varsigma_A}-Z_{\varsigma_A}$ counts the progeny
of $Z_{\varsigma_A}$ independent particles, and thus with a large
probability $S_{\varsigma_A}$ is not very different from
$Z_{\varsigma_A} \bigl(1+E_\omega
(Y_{\varsigma_A})\bigr)=Z_{\varsigma_A}
R\bigl(\theta^{-\varsigma_A}\omega\bigr),$ where the random
variable $R$ is defined by \eqref{sf}. We will obtain \beq \lim_{t
\to \infty} t^\kappa \pp_+ \bigl( \witi W \geq t \bigr) &=&\lim_{A
\to \infty} \lim_{t\to \infty} t^\kappa \pp_+
\bigl(S_{\varsigma_A} \geq t,\varsigma_A<\tilde \nu \bigr)=\lim_{
A\to \infty} \ee_+ \bigl(Z_{\varsigma_A}^\kappa
K(\theta^{-\varsigma_A} \omega); \varsigma_A<\tilde \nu
\bigr),\feq where the random variable $K(\omega)$ is defined by
\eqref{gvul}.

We shall then end the proof by showing that for all $t$ and $A$
large enough, $\ee_+ \bigl(Z_{\varsigma_A}^\kappa ;
\varsigma_A<\tilde \nu \bigr)$ and therefore $t^\kappa \pp_+
\bigl( \witi W \geq t \bigr) \approx \ee_+
\bigl(Z_{\varsigma_A}^\kappa K(\theta^{-\varsigma_A} \omega);
\varsigma_A<\tilde \nu \bigr)$ is uniformly bounded away from zero
and infinity by constants independent of $\omega.$

To carry out this outline, the three terms in the right-hand side
of \eqref{piruk-matrif} are evaluated in the following series of
lemmas, which are versions of the corresponding statements (Lemmas
2--5) in \cite{kks}, and their proofs are deferred to the end of
this Appendix.

We start with the following corollary to Proposition \ref{moment}.
\begin{lemma} \label{uest3} Assume that Condition B is satisfied.
Then, \item [(a)] There exist $C_3,~C_4>0$ such that $P-\as,$
$\pp( \tilde \nu>n |\calf_0) \leq C_3 e^{-C_4 n},$ for any $n>0.$
\item [(b)] There exists a deterministic function $\eta_t>0,~t
\geq 0$ such that $\lim_{t \to \infty} \eta_t=0$ and
$\pp_+\bigl(\varsigma_A<\tilde \nu\bigr) \leq \eta_A. $
\end{lemma}
Fix now any $\delta>0.$ It follows from part {\em (a)} of Lemma
\ref{uest3} that for any $A>0,$ \beq
\pp_+\left(\sum_{n=0}^{\min\{\varsigma_A,\tilde \nu\}-1} Z_n \geq
\delta t  \right) \leq \pp_+(A\tilde \nu \geq \delta t)\leq C_3 e^{-C_4
\delta t/A}=o(t^{-\kappa}),~~t \ra \infty, \feq and thus \beqn
\label{dereh1} \pp_+(\witi W \geq \delta t,\varsigma_A \geq \tilde
\nu) \leq \pp_+(A\tilde \nu \geq \delta t)\leq C_3 e^{-C_4 \delta
t/A}=o(t^{-\kappa}),~~t \ra \infty, \feqn \beqn \label{dereh2}
\pp_+\left(\suml_{n=0}^{\varsigma_A-1}Z_n \geq \delta
t,\varsigma_A < \tilde \nu\right) \leq \pp_+(A\tilde \nu \geq
\delta t) \leq C_3 e^{-C_4\delta t/A}=o(t^{-\kappa}),~~t \ra
\infty. \feqn
\begin{lemma}
\label{vai} \item [(i)] There exists a constant $K_9>0$ such that
$\pp_+\left(Y_0 \geq t \right) \leq K_9 t^{-\kappa}$ for all
$t>0.$ \item[(ii)] For all $\delta>0$ there exists an
$A_0=A_0(\delta)<\infty$ such that \beqn \label{cor17}
\pp_+\left(\sum_{\varsigma_A \leq n <\tilde \nu }Y_n \geq \delta
t\right) \leq \delta t^{-\kappa}~~~~~\mbox{for all}~A \geq A_0.
\feqn
\end{lemma}
It follows from \eqref{piruk-matrif}, taking estimates
\eqref{dereh1}, \eqref{dereh2} and \eqref{cor17} into account,
that for any $A>A_0(\delta)$ (where $A_0$ is given by
\eqref{cor17}) there exists $t_A>0$ such that \beqn
\label{approx-ww} \pp_+(\varsigma_A <\tilde \nu,S_{\varsigma_A}
\geq t) \leq \pp_+(\witi W \geq t) \leq \pp_+(\varsigma_A<\tilde
\nu,S_{\varsigma_A} \geq t(1-2\delta))+3\delta t^{-\kappa}, \feqn
for all $t>t_A.$ Thus, $\witi W$ can be approximated by
$S_{\varsigma_A}.$

Recall the random variable $R$ defined by \eqref{sf}. Note that
$R(\omega)=E_\omega(Y_0),$ and, denote (as in \eqref{kohav})
$R^{\varsigma_A}= R(\theta^{-\varsigma_A}\omega).$ We have the
following law of large numbers with random normalizing constant
$Z_{\varsigma_A}.$
\begin{lemma}
\label{thelast} \item [(i)] There exist functions
$K_{10}=K_{10}(A)>0$ and $K_{11}=K_{11}(A)>0$ independent of
$\omega$ such that \beqn \label{uest1}
K_{10}(A)<\ee_+\bigl(Z_{\varsigma_A}^\kappa; \varsigma_A<\tilde
\nu\bigr) \leq K_{11}(A). \feqn \item [(ii)] For all $\delta>0$
there exists an $A_1=A_1(\delta)$ such that \beqn \label{new1}
\pp_+ \left(\bigl|S_{\varsigma_A} -Z_{\varsigma_A}
R^{\varsigma_A}\bigr| \geq \delta t, \varsigma_A
  <\tilde \nu\right) \leq \delta t^{-\kappa}
\ee_+\bigl(Z_{\varsigma_A}^\kappa; \varsigma_A<\tilde
\nu\bigr)\feqn for $A \geq A_1.$
\end{lemma}
It follows from \eqref{approx-ww} and \eqref{new1} that for $A$
and $t$ sufficiently large, \beqn \label{approx-w}
&&\pp_+\bigl(\varsigma_A <\tilde \nu, Z_{\varsigma_A}
R^{\varsigma_A}\geq (1+\delta)t\bigr)-\delta
t^{-\kappa}\ee_+\bigl(Z_{\varsigma_A}^\kappa;\varsigma_A<\tilde
\nu\bigr) \leq \pp_+\bigl(\witi W \geq t\bigr)  \nonumber \\&&
\leq \pp_+\bigl(\varsigma_A<\tilde \nu, Z_{\varsigma_A}
R^{\varsigma_A} \geq (1-3\delta)t\bigr)+ \delta
t^{-\kappa}\left(3+
\ee_+\bigl(Z_{\varsigma_A}^\kappa;\varsigma_A<\tilde
\nu\bigr)\right).\feqn For a fixed $A>0,$
we obtain from Condition $\mbox{C}_\kappa$ and the dominated
convergence theorem that \beqn \label{1700} &&\lim_{t \ra
\infty}t^{\kappa}\pp_+\bigl(\varsigma_A<\tilde \nu,Z_{\varsigma_A}
R^{\varsigma_A} \geq t\bigr)= \lim_{t \ra \infty} t^{\kappa}
\ee_+\left(I(\varsigma_A<\tilde \nu) \cdot \pp_+\bigl(
Z_{\varsigma_A} R^{\varsigma_A} \geq t \bigl|\calf_{\varsigma_A}
\bigr)\right) \nonumber\\&&= \ee_+\bigr(I(\varsigma_A<\tilde
\nu)\cdot Z_{\varsigma_A}^\kappa \cdot
K\bigl(\theta^{-\varsigma_A}
\omega\bigr)\bigr)=\ee_+\bigr(Z_{\varsigma_A}^\kappa \cdot
K\bigl(\theta^{-\varsigma_A} \omega\bigr);\varsigma_A<\tilde
\nu\bigr), \feqn and, with constants $K_3$ and $K_4$ defined in
\eqref{bound}, \beq
K_3\ee_+\bigr(Z_{\varsigma_A}^\kappa;\varsigma_A<\tilde \nu\bigr)
\leq t^{\kappa}\pp_+\bigl(\varsigma_A<\tilde \nu,Z_{\varsigma_A}
R^{\varsigma_A} \geq t\bigr) \leq
K_4\ee_+\bigr(Z_{\varsigma_A}^\kappa ;\varsigma_A<\tilde
\nu\bigr)\feq for all $t$ sufficiently large.

It follows from \eqref{approx-w} and \eqref{1700} that \beq
\lim_{t \ra \infty} t^{\kappa}\pp_+\bigl(\witi W \geq t\bigr)
=\lim_{A \ra \infty} \ee_+\bigr(Z_{\varsigma_A}^\kappa \cdot
K\bigl(\theta^{-\varsigma_A} \omega\bigr) ;\varsigma_A<\tilde
\nu\bigr), \feq where the last limit is finite by \eqref{bound}
and \eqref{uest1}. The limit in the right-hand side exists since
the limit in the left-hand side does not depend of $A.$

Furthermore, it follows from \eqref{approx-w} and \eqref{bound}
that for some $\delta_0>0,A_2>0,$ \beq
0<\left(\fracd{K_3}{(1+\delta_0)^\kappa}-\delta_0\right)\cdot
\ee_+\bigl(Z_{\varsigma_A}^\kappa;\varsigma_A<\tilde \nu\bigr)
&\leq& t^\kappa \pp_+\bigl(\witi W \geq t\bigr) \\&\leq&
\left(\fracd{K_4}{(1-3\delta_0)^\kappa}+\delta_0\right)\cdot
\ee_+\bigl(Z_{\varsigma_A}^\kappa;\varsigma_A<\tilde
\nu\bigr)+3\delta_0 ,\feq for all $t>t_0.$ Therefore, by
$\eqref{uest1},$ \beq
0<K_{10}(A_2)\left(\fracd{K_3}{(1+\delta_0)^\kappa}-\delta_0
\right) \leq t^\kappa \pp_+\bigl(\witi W \geq t\bigr) \leq
K_{11}(A_2)\left(\fracd{K_4}{(1-3\delta_0)^\kappa}+\delta_0\right)+
3\delta_0,\feq completing the proof of Proposition \ref{st39}.
\subsection*{Proof of Lemma \ref{uest3}}
\label{proof-uest3}
{\em (a)} Follows from part {\em (a)} of Lemma \ref{st17} (which
itself is a corollary to Proposition \ref{moment}).\\ \centerline{}\\
{\em (b)} It is enough to consider $A \in \nn.$ For any $n>0$ we
have \beqn \label{uniset} \pp_+(\varsigma_A <\tilde
\nu)&=&\pp_+(\varsigma_A <\tilde \nu,\tilde \nu
>n)+\pp_+(\varsigma_A <\tilde \nu, \tilde \nu \leq n)
\leq \pp_+(\tilde \nu >n)+\pp_+(\varsigma_A < n)  \nonumber \\
&\leq& C_3e^{-C_4n}+\pp_+(\varsigma_A < n).\feqn For any $n \in
\nn$ let $b_n =(1-1/n)^{1/n}$ and define a sequence of natural
numbers $\{a_{i,n}\}_{i=0}^n$ by the following rule: $a_{0,n}=0$
and \beq a_{i+1,n}=\min\left\{j \in \nn:j
>\max\left\{a_{n-1,n-1};\fracd{(a_{i,n}+1)(1-\epsilon)} {(1-b_n)
\epsilon}\right\} \right\}.\feq  Then, \beq &&\pp_+\bigl(Z_i
>a_{i,n} |Z_j \leq a_{j,n},~ j=0,1,\ldots, i-1\bigr)\leq
\fracd{1}{a_{i,n}} \ee_+\bigl(Z_i |Z_{i-1}=a_{i-1,n}\bigr)\\&&=
\fracd{(a_{i-1,n}+1)(1-\omega_{-i+1})} {a_{i,n} \cdot
\omega_{-i+1}} \leq \fracd{(a_{i-1,n}+1)(1-\epsilon)} {a_{i,n}
\cdot \epsilon} \leq 1-b_n. \feq We conclude that \beq
\pp_+\bigl(Z_i \leq a_{i,n} |Z_j \leq a_{j,n},~ j=0,1,\ldots,
i-1\bigr) \geq b_n, \feq and hence $\pp_+\bigl(\varsigma_A
(a_{n,n})>n\bigr) \geq \pp_+\bigl(Z_i \leq
a_{i,n},~i=1,2,\ldots,n\bigr) \geq 1-1/n.$  By construction,
$a_{n,n}$ is a strictly increasing sequence and it follows from
\eqref{uniset} that for any $A>a_{n,n},$ \beq
\pp_+\bigl(\varsigma_A (A)<\tilde \nu \bigr)\leq
\pp_+\bigl(\varsigma_A (a_{n,n})<\tilde \nu \bigr) \leq C_3
e^{-C_4n}+1/n,\feq completing the proof. \qed
\subsection*{Proof of Lemma \ref{vai}}
\label{proof-vai} {\em (i)} Recall $R^n=1+\sum\limits_{i=1}^\infty
\prod\limits_{j=n}^{n+i-1} \rho_{-j}$ and let $A_n=
Z_{0,n}-Z_{0,n-1} \rho_{-(n-1)}.$ Then, $Y_0=
\sum\limits_{n=1}^\infty A_n R^n ,$ and using the identity
$\sum_{n=1}^\infty n^{-2}=\pi^2/6 <2,$ we obtain from Condition
$\mbox{C}_\kappa$ that \beq \pp_+(Y_0 \geq
t)&=&\pp_+\left(\sum_{n=1}^\infty A_n R^n \geq 6 \pi^{-2}
t\sum_{n=1}^\infty n^{-2}\right) \leq \sum_{n=1}^\infty
\pp_+\left(|A_n|R^n \geq \frac{t}{2n^2} \right)  \\&\leq& 2^\kappa
t^{-\kappa} K_4 \sum_{n=1}^\infty n^{2\kappa}
\ee_+\bigl(|A_n|^\kappa \bigr). \feq Since (cf. \cite[pp.
158--159]{kks})~ $\ee_+\bigl(|A_n|^\kappa \bigr) \leq K_{12} E_P\left(
\prod_{i=0}^{n-2}\rho_{-i}^{\kappa/2}\big|
\calf_0\right)$ for some constant $K_{12}>0,$ it follows from
Condition B that
$\pp_+(Y_0 \geq t) \leq K_9 t^{-\kappa},$ for some $K_9>0.$ \\
\centerline{} \\
\noindent {\em (ii)} Recall the $\sigma$-algebra $\calf_n$ defined
in \eqref{sigmag}. Using the first part of the proposition, we
obtain: \beq &&\pp_+\left(\sum_{\varsigma_A \leq n <\tilde \nu
}Y_n \geq \delta t \right) = \pp_+\left(\sum\limits_{n=1}^\infty
Y_n I(\varsigma_A \leq n <\tilde \nu) \geq 6\delta t \pi^{-2}
\sum\limits_{n=1}^\infty n^{-2}\right) \\ &&\leq
\sum\limits_{n=1}^\infty \ee_+ \Bigl(I(\varsigma_A \leq n <\tilde
\nu)\cdot \pp \bigl(Y_n \geq 1/2 \cdot \delta t n^{-2}
\bigl|\calf_n\bigr) \Bigr)  \\ &&\leq K_9 2^\kappa t^{-\kappa}
\delta^{-\kappa}\ee_+\left(\tilde \nu^{2\kappa+1}; \varsigma_A
<\tilde \nu\right) \leq K_9 2^\kappa t^{-\kappa}
\delta^{-\kappa}\sqrt{\ee_+\left(\tilde
\nu^{4\kappa+2}\right)}\cdot \sqrt{\pp_+(\varsigma_A <\tilde \nu)}
. \feq The claim follows now from Lemma \ref{uest3}, the first
square root being bounded and the second one going to zero as $A
\ra \infty,$ both uniformly in $\omega.$ \qed
\subsection*{Proof of Lemma \ref{thelast}} \label{proof-thelast}
{\em (i)} For the lower bound: \beq
\ee_+\bigl(Z_{\varsigma_A}^\kappa; \varsigma_A<\tilde \nu\bigr)
&\geq& A^\kappa \pp_+\bigl(\varsigma_A<\tilde \nu\bigr) \geq
A^\kappa \pp_+\bigl(Z_1=A+1\bigr) =A^\kappa
\omega_0(1-\omega_0)^{1+A}
\\&\geq& A^\kappa \epsilon^{A+2}:=K_5(A)>0. \feq We now turn to
the upper bound. For a fixed environment $\omega$ we obtain, by
using the Markov property of $Z_n$ in the second equality and the
ellipticity condition {\em (B1)} in the last two inequalities,
\beq E_\omega \left( Z_{\varsigma_A}^\kappa \right) &=&\sum_{n
\geq 1} \sum_{a=0}^A E_\omega \left( Z_n^\kappa \left| \right.
\varsigma_A =n, Z_n>A, Z_{n-1}=a\right) P_\omega\left(\varsigma_A=n,
Z_{n-1}=a \right)=
\\&=& \sum_{n \geq 1} \sum_{a=0}^A E_\omega \left( Z_n^\kappa \left| \right.
Z_n>A, Z_{n-1}=a\right)P_\omega\left(\varsigma_A=n,
Z_{n-1}=a \right)\\
&\leq& \sup_{\omega,n \in \nn,a \leq A} E_\omega \left( Z_n^\kappa
\left| \right. Z_n>A,Z_{n-1}=a \right) \leq \sup_{\omega,n \in
\nn,a \leq A}\fracd{E_\omega \left( Z_n^\kappa \left| \right.
Z_{n-1}=a \right) }{P_\omega \left( Z_n >A \left| \right.
Z_{n-1}=a \right) } \\&\leq& \sup_\omega \fracd{E_\omega \left(
Z_1^\kappa \left| \right. Z_0=A \right) }{P_\omega \left( Z_1 >A
\left| \right. Z_0=0 \right) } \leq (A+1)^{\kappa+1}
\epsilon^{-A-2}\sup_\omega E_\omega
\left[\left(V_{0,0}\right)^\kappa \right] <\infty ,\feq where the
random variables $V_{n,j}$ are defined in \eqref{bigv}. This
completes the proof of part {\em(i)} of the Lemma.
\\ \centerline{} \\
\noindent{\em (ii)} The proof is similar to that of Lemma
\ref{vai}. If $\varsigma_A<\tilde \nu,$ let \beq
S_{\varsigma_A,j}=\mbox{number of progeny alive at
time}~j~\mbox{of the $Z_{\varsigma_A}$ particles present at
time}~\varsigma_A,\feq and
$B_j=S_{\varsigma_A,j}-S_{\varsigma_A,j-1} \cdot \rho_{-(j-1)}.$
We have $\sum_{j=\varsigma_A}^\infty S_{\varsigma_A,j}
-Z_{\varsigma_A} R^{\varsigma_A} =\sum_{j=\varsigma_A}^\infty B_j
R^j,$ and obtain from Condition $\mbox{C}_\kappa$ that on the set
$\{\varsigma_A < \tilde \nu\},$ \beq &&\pp_+\left( \Bigl|
\sum_{j=\varsigma_A}^\infty S_{\varsigma_A,j} -Z_{\varsigma_A}
R^{\varsigma_A} \Bigr| \geq \delta t \Bigl| \calf_{\varsigma_A}
\right)\leq \sum_{j=\varsigma_A}^\infty \ee_+ \bigl( \pp_+
\bigl(|B_j| R^j \geq \fracd{\delta
t}{2(j-\varsigma_A+1)^2} \bigl| B_j,\calf_{\varsigma_A} \bigr)\bigr) \\
&&\leq K_4 \left(\fracd{2}{\delta t}\right)^\kappa
\sum_{n=0}^\infty (n+1)^2 \cdot \ee_+ \bigl( \bigl|
B_{\varsigma_A+n}\bigr|^\kappa~
\bigl|\calf_{\varsigma_A}\bigr).\feq Since (cf. \cite[p. 164]{kks})~
$\ee_+ \bigl( \bigl| B_{\varsigma_A+n}\bigr|^\kappa~
\bigl|\calf_{\varsigma_A}\bigr) \leq K_{13}Z_{\varsigma_A}^\kappa
E_P \left(\prod_{i=\varsigma_A}^{j-2}
\rho_{-i}^{\kappa/2}\right),$ it follows from Condition
$\mbox{C}_\kappa$ that for some $K_{14}>0,$ \beq &&\pp_+
\left(\Bigl| \sum_{j=\varsigma_A}^\infty S_{\varsigma_A,j}
-Z_{\varsigma_A} \witi S_{\varsigma_A}\Bigr| \geq \delta t;
\varsigma_A <\tilde \nu\right) \leq \left(\fracd{K_{14}}{t\delta}
\right)^\kappa \ee_+
\bigl(Z_{\varsigma_A}^{\kappa/2};\varsigma_A<\tilde
\nu\bigr)\\&&\leq \left(\fracd{K_{14}}{t\delta \sqrt{A}}
\right)^\kappa\ee_+
\bigl(Z_{\varsigma_A}^\kappa;\varsigma_A<\tilde \nu\bigr)\leq
\delta t^{-\kappa}\ee_+
\bigl(Z_{\varsigma_A}^\kappa;\varsigma_A<\tilde \nu \bigr), \feq
for $A \geq A_2(\delta).$ \qed
\end{appendix}
\section*{Acknowledgment}
We wish to thank Julien Bremont for his careful reading of a
preliminary version of the paper and for several valuable suggestions.
\providecommand{\bysame}{\leavevmode\hbox
to3em{\hrulefill}\thinspace}
\providecommand{\MR}{\relax\ifhmode\unskip\space\fi MR }
\providecommand{\MRhref}[2]{%
  \href{http://www.ams.org/mathscinet-getitem?mr=#1}{#2}
} \providecommand{\href}[2]{#2}

 \end{document}